\documentclass[10pt,letter]{article}

\usepackage{placeins}
\usepackage{multicol}

\RequirePackage{hyperref}
\RequirePackage[numbers]{natbib}

\usepackage{enumitem}
\usepackage{titlesec}

\usepackage{amsmath}
\usepackage{amsthm}
\usepackage{amssymb}

\usepackage{bbm}
\usepackage{bm}
\usepackage{nicefrac}

\usepackage{mathtools}

\usepackage{tikz}

\usepackage{epsfig}
\usepackage{graphicx}
\usepackage[export]{adjustbox}

\usepackage{resizegather}
\usepackage{color}

\usepackage{float}
\usepackage{wrapfig}
\usepackage{caption}
\usepackage{subcaption}

\usepackage{comment}
 \usepackage{appendix}

\usepackage{algorithm}
\usepackage{algorithmic}
\usepackage{listings}

\newtheorem{theorem}{Theorem}
\newtheorem{corollary}{Corollary}

\newtheorem*{theorem*}{Theorem}
\newtheorem*{lemma*}{Lemma}
\newtheorem*{problem*}{Problem}

\newenvironment{remark}[1][Remark]{\begin{trivlist}
\item[\hskip \labelsep {\bfseries #1}]}{\end{trivlist}}

\graphicspath{{./figures/}}

\DeclareMathOperator{\Ii}{\mathcal{I}}

\DeclareMathOperator{\EE}{\mathbb{E}}
\DeclareMathOperator{\RR}{\mathbb{R}}

\DeclareMathOperator{\Prob}{\mathbb{P}}

\DeclareMathOperator{\uu}{\mathbf{u}}
\DeclareMathOperator{\vv}{\mathbf{v}}

\DeclareMathOperator{\yy}{\mathbf{y}}
\DeclareMathOperator{\zz}{\mathbf{z}}

\DeclareMathOperator{\FF}{\mathbf{F}}
\DeclareMathOperator{\GG}{\mathbf{G}}
\DeclareMathOperator{\XX}{\mathbf{X}}

\DeclareMathOperator{\bmu}{\bm{\mu}}
\DeclareMathOperator{\bSig}{\bm{\Sigma}}

\usetikzlibrary{shapes, arrows, positioning}

\tikzstyle{decision} = [diamond, draw, fill = green!20, 
    text width = 4.5em, text badly centered, node distance = 3cm, inner sep = 0pt]
\tikzstyle{block} = [rectangle, draw, fill = blue!20, 
    text width = 10em, text centered, rounded corners, minimum height = 4em]
\tikzstyle{line} = [draw, -latex']
\tikzstyle{cloud} = [draw, ellipse, fill = red!20, node distance = 3cm,
    minimum height = 2em]
\tikzstyle{circ} = [circle, draw, fill = red!20, 
    text width = 10em, text centered, minimum height = 2em]

\begin{document}

\title{Sparse Equisigned PCA:\\ Algorithms and Performance Bounds in the Noisy Rank-1 Setting\footnote{This work was partially supported by the NSF and the ONR.}}
\author{Arvind Prasadan\footnote{University of Michigan, Dept of EECS, prasadan{@}umich.edu}, Raj Rao Nadakuditi \footnote{University of Michigan, Dept of EECS, rajnrao{@}umich.edu}, and Debashis Paul\footnote{University of California, Davis, Department of Statistics, debpaul{@}ucdavis.edu}}

\maketitle

\begin{abstract}
Singular value decomposition (SVD) based principal component analysis (PCA) breaks down in the high-dimensional and limited sample size regime below a certain critical eigen-SNR that depends  on the dimensionality of the system and the number of samples. Below this critical eigen-SNR, the estimates returned by the SVD are asymptotically uncorrelated with the latent principal components. We consider a setting where the left singular vector of the underlying rank one signal matrix is assumed to be sparse and the right singular vector is assumed to be equisigned, that is, having either only nonnegative or only nonpositive entries. We consider six different algorithms for estimating the sparse principal component based on different statistical criteria and prove that by exploiting sparsity, we recover consistent estimates in the low eigen-SNR regime where the SVD fails. Our analysis reveals conditions under which a coordinate selection scheme based on a \textit{sum-type decision statistic} outperforms schemes that utilize the $\ell_1$ and $\ell_2$ norm-based statistics. We derive lower bounds on the size of detectable coordinates of the principal left singular vector and utilize these lower bounds to derive lower bounds on the worst-case risk. Finally, we verify our findings with numerical simulations and a illustrate the performance with a video data where the interest is in identifying objects.
\end{abstract}

\section{Introduction} \label{sec:intro}

It is well-understood that singular value decomposition (SVD) based principal component analysis (PCA) breaks down in the high-dimensional and limited sample size regime below a certain critical eigen-SNR (eigenvalue signal-to-noise ratio) that depends on the dimensionality of the system and the number of samples \cite{johnstone2009consistency, birnbaum2013minimax}. Several sparse PCA algorithms have been proposed in the literature (see \cite{johnstone2009consistency, birnbaum2013minimax, d2007direct, ma2013sparse, yuan2013truncated, berthet2013optimal}) and have been shown to successfully estimate the principal components in the low eigen-SNR regime where the SVD fails. 

Prior work in this area primarily considers the 
Gaussian signal-plus-noise model with
random effects, where the signal matrix is assumed to have sparse left singular vectors, normally distributed right singular vectors, and the noise matrix 
is assumed to have normally distributed \emph{i.i.d.} entries. Here, we consider the setting where the left singular vector of the rank one signal matrix is sparse \emph{and} the right singular vector is assumed to be equisigned. We say that a vector is equisigned if its entries are all non-negative or all non-positive. This is motivated by applications such as diffusion imaging in MRI where the right singular vector represents a physical quantity (e.g. intensity as the diffusion agent is absorbed by a tissue) that is non-negative, by imaging problems such as foreground-background separation in video data \cite{piccardi2004background, vaswani2018robust} and object detection in astronomy \cite{ren2018non}, where the data are naturally non-negative, and by problems in bioinformatics where the data are (non-negative) counts of genes \cite{taslaman2012framework}. When analyzing data that are non-negative, it is logical to take advantage of this property, and investigate how we may use this knowledge to do better than the (generic) alternatives. Alternatively, a practitioner may seek to use techniques that constrain or impose non-negativity to preserve interpretability of the results, e.g., non-negative matrix factorization. Additionally, we motivate the rank-$1$ assumption by noting that for a video with a static background, the foreground is a perturbation of a rank-$1$ background \cite{moore2019panoramic, gao2017augmented}. Finally, even though we do not pursue this angle here, our framework can be extended to deal with the scenario where
the signal can be viewed of a rank 1 tensor with all but one of the representors in the Kroneker product representation of the tensor is an equisigned vector.

There is precedent for and prior work on non-negative PCA, including the sparse biased PCA in \cite{brennan2018reducibility}, the sparse PCA with non-negativity priors in \cite{perry2018optimality}, and the work in \cite{montanari2015non}. These works differ from our work in that they impose non-negativity on the factors or the left singular vectors. In this work, we study sparse factors with non-negative loadings; i.e., we are solving a different problem in this work.

A natural question at this juncture is the following: \emph{how does our problem differ from that solved by Non-Negative Matrix Factorization (NNMF)?} NNMF takes a given matrix $\XX$ and looks for non-negative matrices $\FF$ and $\GG$ such that $\XX = \FF\GG^T$ \cite{huang2014non, wang2013nonnegative}. Ordinary NNMF has no sparsity constraints. We might impose such constraints, as is done in \cite{hoyer2004non} and \cite{liu2012constrained}, but except in special cases, these solutions have no known theoretical guarantee of statistical performance.
This problem partly stems from the fact that solutions to the corresponding optimization problems may not be unique. In contrast, our problem only constrains the right singular vectors, while the left singular vectors are free to take any sign. The work in \cite{ding2010convex} extends the NNMF framework to one wherein only one of the factors is non-negative; nevertheless, the rest of the constraints we impose are not included. The work in \cite{zass2007nonnegative} seeks factors (left singular vectors) with disjoint supports and non-negative loadings, but this definition of sparsity does not match that from the sparse PCA literature. Hence, NNMF is not an answer to the problem we consider herein. 

The main contribution of this paper is a rigorous sparsistency analysis of the various algorithms that brings into focus the various very-low eigen-SNR regimes where the new algorithms work and the SVD based methods provably fail. Additionally, a major novelty of this work is the integration of FDR-controlling (False Discovery Rate) hypothesis testing to the Sparse PCA problem.

Our analysis illustrates the situations where the sum based coordinate selection scheme dramatically outperforms the $\ell_1$ and $\ell_2$ \cite{johnstone2009consistency, birnbaum2013minimax} based sparse PCA schemes. Additionally, our proposed algorithms are non-iterative, do not require the computation of the sample covariance matrix, and do not require knowledge 
of the sparsity level. We separate our algorithms into two groups: one where the Family-Wise Error Rate (FWER) is controlled, and another where the False Discovery Rate (FDR) is controlled. We utilize sharp tail probability bounds for relevent statistics to derive our FWER-controlling estimators \cite{boucheron2013concentration}. For the FDR controlling estimators, we relate the problem at hand to that of the sparse normal means problem \cite{donoho2004higher}.

This paper is organized as follows. In Section \ref{sec:algorithms}, we describe three algorithms for estimating the sparse principal component  that utilize a coordinate selection scheme based on the sum, $\ell_1$, and $\ell_2$ norm-based statistics respectively. We call our family of algorithms \emph{SEPCA}, an abbreviation for Sparse Equisigned PCA. Section \ref{sec:FDR} proposes three FDR-controlling refinements of the sum- and $\ell_2$-based algorithms in Section \ref{sec:algorithms} by relating coordinate detection to the sparse normal means estimation problem.  In Section \ref{sec:risk} we show how the estimation performance is governed by the size of the smallest detectable coordinate, which we analyze in Section \ref{sec:perf} and validate using numerical simulations in Section \ref{sec:sim}. In Section \ref{sec:comp}, we provide some geometric intuitions about the relative performance of three of our algorithms. We show that the sum statistic is potentially the most powerful, while the $\ell_1$ is the least powerful. We provide some concluding remarks in Section \ref{sec:conclusions}. 

\section{Problem Formulation} \label{sec:model}

Let $\XX \in \mathbb{R}^{p \times n}$ be a real-valued signal-plus-noise data matrix of the form
\begin{equation} \label{eqn:model}
\XX = \theta \uu \vv^T + \sigma \GG.
\end{equation}
The columns of the $p \times n$ data matrix $\XX$ represent $p$-dimensional observations. In (\ref{eqn:model}), $\uu$ and $\vv$ are the left and right singular vectors of the rank-one latent signal matrix, and have entries $u_i$ and $v_j$, respectively. The entries of $\GG$, the noise matrix, are assumed to be \emph{i.i.d.} Gaussian random variables with mean $0$ and variance ${1}/{n}$. We assume that $\uu \in \mathbb{R}^p$ has unit norm and is sparse in the sense of small $\ell_0$ norm, with $s \ll p$ non-zero entries, where $s / n \rightarrow 0$. That is, for a set $I = \lbrace i_1, \cdots, i_s \rbrace \subset \lbrace 1, \cdots, p\rbrace$, 
\begin{equation} \label{eqn:u}
\begin{array}{ll}
u_i \neq 0  & \text{for $i \in I$,} \\
u_i = 0 & \text{for $i \in I^C$,}
\end{array}
\end{equation}
where $I^C$ denotes the complement of $I$. We further assume $\vv \in \mathbb{R}^n$ to be of unit norm, deterministic, and equisigned. Given $\XX$, our goal is to recover $\uu$ and $\vv$.

Note that the $(i, k)$ entry of $\XX$, $X_{ik}$, is a Gaussian random variable with mean $[\theta u_i] v_k$ and variance ${\sigma^2}/{n}$. Moreover, it follows that
$$\mathbb{E}(\XX\XX^T) = \theta^2 \uu \uu^T + {\sigma^2} \Ii_p,$$ 
where $\Ii_p$ denotes the $p \times p$ identity matrix.
The quantity $\left(\theta / \sigma\right)^2$ is, for this model, the eigen-SNR (signal-to-noise ratio). 

\subsection{Motivation: Breakdown of PCA / SVD} \label{ssec:breakdown}

From \cite{FBG_RRN_2012}, we have the following result: let $\widehat{\uu}$ be the estimate of $\uu$ given by the Singular Value Decomposition (SVD) of $\XX$, and let ${p(n)}/{n}$ have limit $c \in [0, \infty]$ as $n$ grows, with $\theta$ fixed and $\sigma = 1$. Then, with probability $1$,
\begin{equation} \label{eqn:break}
\left|\langle \widehat{\uu}, \uu\rangle\right|^2 \rightarrow \left\lbrace \begin{array}{ll} 1 - \frac{c \left(1 + \theta^2\right)}{\theta^2 \left(c + \theta^2\right)} & \text{if $\theta \geq c^{1/4}$,} \\ 0 & \text{otherwise}.\end{array} \right.
\end{equation}
For general $\sigma$, we replace $\theta$ by $\theta/\sigma$ in (\ref{eqn:break}).
Hence, SVD based PCA leads to inconsistent estimates of $\uu$ (and also for $\vv$, which can be deduced from (\ref{eqn:break})) when the dimension $p$ is comparable to or larger 
than the sample size $n$. Moreover, in the low eigen-SNR regime, the estimates break down completely. SVD does not exploit any assumed structure in $\uu$ and $\vv$. Consequently, (\ref{eqn:break}) holds for arbitrary $\uu$ and $\vv$, including our setting where $\uu$ is sparse and/or $\vv$ is equisigned. Our goal, in what follows, is to derive consistent estimators for $\uu$ and $\vv$ that outperform the SVD by exploiting the sparsity of $\uu$ and the equisigned nature of $\vv$.

\subsection{Problem Statement}

{Note that we have assumed that $\uu$ is sparse and that the sparsity $s$ is such that $s / n$ has limit zero. Hence, if we had oracle knowledge of the sparsity pattern $I$ (the indices of $\uu$ that have non-zero coordinates), restricting the matrix $X$ to those rows indexed by $I$ and performing the SVD on the smaller matrix would yield a consistent estimator for the non-zero elements of $\uu$ and the vector $\vv$. This conclusion follows from (\ref{eqn:break}), since the value $c$ is replaced with $s / n$, which has limit zero. Thus, if we derived consistent estimators of the support of $\uu$, we have a consistent two-stage estimation procedure of the vectors $\uu$ and $\vv$.}

{Formally, we are interested in finding a procedure that estimates $I$ by $\widehat{I}$ such that in the limit $n \rightarrow \infty$,
\begin{equation}
    \left\{\begin{array}{ll}
        \Prob\left(i \in \widehat{I}\right) \rightarrow 1 & \textrm{ if $i \in I$}, \\
        \Prob\left(i \in \widehat{I}\right) \rightarrow 0 & \textrm{ if $i \notin I$}.  
    \end{array}\right.
\end{equation}
Equivalently, noting that the Hamming distance of $I$ and $\widehat{I}$, denoted by $d_H\left(I, \widehat{I}\right)$, is given by the cardinality of their symmetric set difference,
$$d_H\left(I, \widehat{I}\right) = \left|\left(I \cup \widehat{I}\right) \backslash \left(I \cap \widehat{I}\right)\right|,$$
we want the expected Hamming distance $\EE d_H\left(I, \widehat{I}\right)$ to have limit $0$, which is stronger
than requiring consistency in recovering the support
(or sparsity pattern) of $\uu$. However, as the work in \cite{butucea2018variable, reeves2008sampling} indicates, this limit will not in general be zero, and will depend on the noise level, signal strength, and sparsity.} 

\section{Proposed Algorithms}\label{sec:algorithms}

We propose six different two-stage algorithms for estimating $\uu$. The first three algorithms 
are designed to control the family-wise error rate (FWER), or, the probability of obtaining a false positive in the coordinate selection. The last three algorithms aim to control the false discovery rate (FDR), or, the proportion of false discoveries (coordinate detections) among all discoveries. We defer discussion of the FDR-based algorithms to Section \ref{sec:FDR}. 

All of the algorithms have the same basic form given in Algorithm \ref{alg:gen}. Given $\XX$, we associate a test statistic $T_i$ to each row of $\XX$. The sparsity of $\uu$ implies that the majority of the rows of $\XX$ are purely noise, so that the majority of the $T_i$ come from the null, noise-only distribution. Hence, based on the statistics $\{T_i\}$, we perform a form of multiple hypotheses testing procedure, and select the set $\widehat{I}$ of indices that are non-null. In this way, we can estimate the \emph{support} of $\uu$, thereby isolating the the rows of $\XX$ that contain the signal. Then, taking the SVD of this submatrix (comprised of only the selected rows of $\XX$) yields a better estimate of the non-zero coordinates in $\uu$, as well as $\vv$. 
\noindent
\begin{algorithm}
\caption{Variable Selection and Estimation Algorithm}
\begin{algorithmic}
\REQUIRE Threshold $\tau_{n, p}$ and form of Test Statistic $T_i$ from Table \ref{tab:details}
\STATE Let $\widehat{I}$ be an empty list
\FORALL{Rows $i$ of $\XX$, $1 \leq i \leq p$} \STATE Form test statistic $T_i$ from row $i$ of $\XX$
\IF {$T_i \geq \tau_{n, p}$} \STATE Add $i$ to $\widehat{I}$ \ENDIF \ENDFOR
\STATE Let $[\widetilde{\uu}, \widetilde{\theta}, \widehat{\vv}] = \text{SVD}(\XX_{\widehat{I}, :})$ be the rank-$1$ SVD of $\XX$ restricted to rows in $\widehat{I} = [i_1, \cdots, i_{|\widehat{I}|}]$
\STATE For $i_k \in \widehat{I}$, let $\widehat{u}_{i_k} =  \widetilde{u}_{k}$; the other entries of $\widehat{\uu}$ are set to $0$.
\end{algorithmic}
\label{alg:gen}
\end{algorithm}

We begin by discussing the FWER-controlling algorithms. The work in \cite{johnstone2009consistency} proposed a covariance thresholding method for Sparse PCA called DT-SPCA; this is equivalent to a coordinate selection scheme based on the $\ell_2$ norm-based statistic. In our terminology and with our choice of thresholds, we label it as \textit{$\ell_2$-SEPCA}. We label the coordinate selection scheme based on the $\ell_1$ norm-based statistic \textit{$\ell_1$-SEPCA}. Finally, the \textit{sum-SEPCA} algorithm utilizes row sums of the data matrix. 

\subsection{Computational Complexity}

{Note that the variable selection part of our procedures has a computational complexity that is $O\left(p n\right)$: the formation of the test statistic is linear in the number of columns, and the formation is repeated once per row. Noting that for a $p \times n$ matrix, the complexity of the rank-$1$ SVD is $O\left(1 \times p n\right)$, we find that if $\left|\widehat{I}\right|$ coordinates are selected, we have an overall complexity of $O\left(p n + \left|\widehat{I}\right| n\right) = O\left(p n\right)$ \cite{allen2016lazysvd}.}

{Computation of the covariance matrix has a (naive) complexity of $O\left(p^2 n\right)$, and in practice is somewhere between $O\left(p^2\right)$ and $O\left(p^3\right)$ \cite{gall2018improved}. Immediately, our methods here are faster than those requiring explicit formation of the covariance matrix \cite{birnbaum2013minimax, johnstone2009consistency, ma2013sparse, yuan2013truncated, berthet2013optimal}. Additionally, there is no iteration or convergence of any optimization problems required. Note that a semi-definite programming-based formulation is at least polynomial in the problem size: $O\left(p^4\right)$ \cite{d2007direct} or $O\left(p^3\right)$ \cite{berthet2013optimal}. The ITSPCA method applied to our rank-$1$ setting would have a cost of $O\left(p s\right)$ per iteration \cite[Sec.~4]{ma2013sparse}. TPower has a similar complexity of $O\left(s p + p\right)$ per iteration \cite{yuan2013truncated}.} 

\subsection{The DT-SPCA Algorithm and Two-Stage Procedures}

{The DT-SPCA algorithm was proposed in \cite{johnstone2009consistency} and later used as the first stage of the ASPCA algorithm given in \cite{birnbaum2013minimax}. The algorithm thresholds the diagonals of the matrix $X X^T$ to perform variable selection: note that in our setting, these values are 
$${\theta^2} u_i^2 + {\sigma^2} \sum_{j = 1}^n G_{ij}^2,$$
with expectation $\left(\theta^2 u_i^2 + \sigma^2\right)$. The DT-SPCA algorithm thresholds these diagonal values at $\sigma^2 \left(1 + \gamma \sqrt{\frac{\log p}{n}}\right)$, where $\gamma > 0$, and then performs PCA on the reduced matrix formed from the selected variables. Noting that the diagonals of $X X^T$ are the same as the row sum-of-squares of $X$, we see that $\ell_2$-SEPCA is essentially the same (up to choice of threshold) as DT-SPCA. }

{However, the innovation of \cite{johnstone2009consistency, birnbaum2013minimax} that we carry forward is the two-stage procedure. That is, we perform some sort of testing to estimate the support of the sparse singular vector $\uu$, and then perform an SVD on the reduced matrix. As we will see in what follows, there is flexibility in the choice of testing or support estimation method. }

\subsection{Statement of Thresholds}

We shall choose the thresholds $\tau_{n, p}$ for the coordinate selection scheme so that in the noise-only case,
\begin{equation} \label{eqn:falsepos}
\mathbb{P} \left(\max_{1 \leq i \leq p} T_i \geq \tau_{n, p}\right) \leq \frac{1}{e p} \rightarrow 0,
\end{equation}
where $e$ is Euler's number, or the base of the natural logarithm. This choice ensures that the probability of a false positive tends to zero as  $p \to \infty$. That is, the FWER is asymptotically zero and is bounded by $1 / ep$ in the finite-dimensional case. Note that the constraint used to control the FWER is simply that the distribution of the noise is log-concave. In the Gaussian case, we obtain the specific expressions given summarized in Table \ref{tab:details}; however, with knowledge of the moments $\mathbb{E} T_i$ and $\text{Var } T_i$, we can repeat our analysis and find thresholds for the $\ell_1$ and $\ell_2$-SEPCA algorithms with \emph{any} log-concave noise distribution. The thresholds are summarized in Table \ref{tab:details}.

\begin{minipage}{\linewidth}
\centering
\captionof{table}{Test Statistics and Thresholds for Algorithm (\ref{alg:gen})} \label{tab:details}
\begin{tabular}{|c|c|c|}
\hline
Algorithm & Statistic $T_i$ & Threshold $\tau_{n, p}$\\
\hline

\rule{0pt}{4ex}   $\ell_1$-SEPCA & $\frac{1}{\sqrt{n}} \sum_{k = 1}^n |X_{i, k}|$ & $\sigma \left(\sqrt{\frac{2}{\pi}} + C_1 \frac{\log e p}{\sqrt{n}}\right)$ \\

\rule{0pt}{4ex} $\ell_2$-SEPCA & $\sum_{k = 1}^n X_{i, k}^2$ & $\sigma^2 \left(1 + C_2  \frac{\log e p}{\sqrt{n}}\right)$ \\

\rule{0pt}{4ex} sum-SEPCA &
$\frac{1}{\sqrt{n}} |\sum_{k = 1}^n X_{i,k}|$ &
$\sigma C_U \sqrt{\frac{\log p}{n}}$ \\

\hline
\end{tabular}
\caption*{See (\ref{eq:thresh_const}) and (\ref{eq:thresh_const_2}) for definitions of the constants $C_2$, $C_1$, and $C_U$.}
\end{minipage}

\begin{remark}
Note that we impose strong control over the FWER and seek to reject individual null hypotheses, instead of weak control and considering the global null hypothesis as in  \cite{berthet2013optimal}.
\end{remark}

\subsection{FWER Thresholds}

\subsubsection{$\ell_2$- and $\ell_1$-SEPCA}

In the noise-only cases, the statistics for $\ell_2$- and $\ell_1$-SEPCA are distributed as scaled $\chi_n^2$ and sums of half-normal, respectively. Both of these quantities are log-concave random variables, so we may apply the result in \cite{Latala2011} to set the threshold $\tau_{n, p}$ in both cases. 

Defining $K$ to be some absolute constant (we may use $K = e$, as in \cite{bobkov2003convex}), we define the constants 
\begin{equation} \label{eq:thresh_const}
C_2 = \sqrt{2} K \text{ and } C_1 = K \sqrt{\left(1 - {2}/{\pi}\right)}.
\end{equation}

\subsubsection{sum-SEPCA}

From Proposition $4.4$ of \cite{BoucheronThomas2012}, we obtain that the threshold for sum-SEPCA is given by 
\begin{gather} \label{eqn:tau_sepca}
\tau_{n, p} = \frac{\sigma}{\sqrt{n}} \left(\sqrt{2 \log p} + \frac{1}{U(p)} \left(\frac{1}{3} \log{e p}+ \sqrt{\log {e p}}\right) + \delta_p\right).
\end{gather}
In (\ref{eqn:tau_sepca}), we have that
\begin{equation}
U(p) = \sqrt{2} \text{ Erf}^{-1}\left(1 - \frac{1}{p}\right) \text{ and } \delta_p \asymp \frac{\pi^2}{12} \left(\log p\right)^{-3/2},
\end{equation}
where $\text{Erf}$ denotes the \emph{error function}, or alternatively, the cumulative distribution function of a standard Gaussian random variable is given by 
\begin{equation}\label{eqn:erf}
\Phi(x) = \frac{1}{2} \left(1 + \text{Erf}\left(\frac{x}{\sqrt{2}}\right)\right).
\end{equation}
Moreover, $\tau_{n, p} \leq \sigma C_{U} \sqrt{\frac{\log p}{n}}$ for some constant $C_U$. For a fixed value of $p$, choosing
\begin{equation} \label{eq:thresh_const_2}
\kappa_U \geq \frac{\sqrt{2}}{U(p)} \left(3 + \sqrt{\log{p}}\right) > 1 \text{ and } C_U = \sqrt{2} + \frac{\kappa_U}{3 \sqrt{2}}
\end{equation}
is sufficient. The choice of $1/ep$ is the largest bound justified by Proposition 4.4 of \cite{BoucheronThomas2012}, so we have calibrated all of our algorithms to the same constant factor times $1/p$. The thresholds are summarized in Table \ref{tab:details}. 

\subsection{Estimation of the Noise Variance, $\sigma^2$}

{In this work, we assume that the noise variance $\sigma^2$ is known; however, in general, estimation of $\sigma^2$ may not be straightforward \cite{passemier2017estimation}. Recently proposed procedures such as those proposed in \cite{passemier2017estimation, pastor2012robust, socheleau2014testing} could be employed to estimate the noise variance, and we point the interested reader to these references for more theoretical background on the problem. We note that in most applications, including the video example we consider, one can obtain a relatively sparse representation of the object in a multiscale basis such as a wavelet basis \cite[Sec.~7.5]{johnstone2013gaussian}. Under such circumstances, under the assumed additive, isotropic noise model, we can easily obtain a consistent estimate of $\sigma^2$ by utilizing the inherent sparsity of the signal, especially in finer scales. This can be done, for example, by computing the variance of the wavelet coefficients in the finest scale \cite[Sec.~7.5]{johnstone2013gaussian}. One can obtain a more robust estimate by taking the median absolute deviation of the coefficients about their median and then by multiplying its square with a known scale factor (assuming normality) \cite{pastor2012robust, johnstone2013gaussian}. }

\section{Controlling the False Discovery Rate}\label{sec:FDR}

So far, we have controlled the probability of a false alarms when detecting coordinates. However, there are two relevant observations to make. First, under the Gaussian noise, rank-1, and equisigned assumptions, the vector of test statistics $\{T_i\}$ in the sum-SEPCA algorithm looks like a sparse vector plus Gaussian noise (or a vector of $\chi_n^2$-variates with varying non-centralities, in the $\ell_2$-SEPCA algorithm). Secondly, controlling the false discovery rate, that is, the proportion of rejected nulls that are false positives, can lead to increased detection power relative to controlling the false positive rate. We hence look at FDR-controlling tests for the \emph{Sparse Normal Means} problem. 

That is, given a vector of test statistics (as before), we replace the thresholding and selection in Algorithm \ref{alg:gen} with an FDR-controlling selection procedure. We summarize this change in Algorithm \ref{alg:gen_fdr}. There are three procedures we consider. The first two are known as Higher Criticism, and directly extend the sum- and $\ell_2$-SEPCA algorithms \cite{donoho2004higher, donoho2015higher}. The third is a method for detection in the sparse normal means problem that comes out of complexity-penalized estimation theory for linear inverse problems \cite{johnstone2014adaptation}.
\noindent
\begin{algorithm}
\caption{FDR-Controlling Variable Selection and Estimation Algorithm}
\begin{algorithmic}
\REQUIRE Test Statistic $T_i$ from Table \ref{tab:details} and Selection Procedure
\STATE Let $\widehat{I}$ be an empty list
\FORALL{Rows $i$ of $\XX$, $1 \leq i \leq p$} \STATE Form test statistic $T_i$ from row $i$ of $\XX$  \ENDFOR
\STATE Perform an FDR-Controlling selection procedure, and add the selected indices to $\widehat{I}$
\STATE Let $[\widetilde{\uu}, \widetilde{\theta}, \widehat{\vv}] = \text{SVD}(\XX_{\widehat{I}, :})$ be the rank-$1$ SVD of $\XX$ restricted to rows in $\widehat{I} = [i_1, \cdots, i_{|\widehat{I}|}]$
\STATE For $i_k \in \widehat{I}$, let $\widehat{u}_{i_k} =  \widetilde{u}_{k}$; the other entries of $\widehat{\uu}$ are set to $0$.
\end{algorithmic}
\label{alg:gen_fdr}
\end{algorithm}


\subsection{Higher Criticism}\label{ssec:HC}

\subsubsection{Formulation} 

Assume we have $p$ independent tests of the form
\begin{equation} \label{eqn:HC_norm_test}
\begin{array}{ll}
H_{o, i}: & W_i \sim \mathcal{N}\left(0, 1\right),  \\
H_{1, i}: & W_i \sim \mathcal{N}\left(\mu_i, 1\right),
\end{array}
\end{equation}
and assume that at most $p^{1 - \beta}$ of the $p$ hypotheses are truly non-null, for some $\beta \in (1/2, 1)$. Further assume that the non-null means have magnitude 
$$\mu_i = \mu_p = \sqrt{2 r \log p},$$
for $r \in (0, 1)$. Here, the means will correspond to the coordinate size. Note that the expected maximum of $p$ standard Gaussian random variables is upper bounded by $\sqrt{2 \log p}$, with the bound being asymptotically sharp. 

If we let $p_{(1)} \leq p_{(2)} \leq \cdots \leq p_{(p)}$ be the sorted p-values of the individual tests, we may define the Higher Criticism statistic: 
\begin{equation} \label{eqn:HC_stat}
HC_p = \max_{i : 1/p \leq p_{(i)} \leq 1/2} \frac{\sqrt{p} \left(i / p - p_{(i)}\right)}{\sqrt{p_{(i)} (1 - p_{(i)})}}.
\end{equation}
Rejecting the global null hypothesis (that there are no non-null coordinates) when $HC_p > \sqrt{2 \log \log p} (1 + o(1))$ leads to asymptotically full power when $r$ is greater than some decision boundary $\rho$, and that under the global null, 
\begin{equation}
\frac{HC_p}{\sqrt{2 \log \log p}} \rightarrow 1
\end{equation}
in probability as $n,p \rightarrow \infty$. The function $\rho$ depends on the sparsity index $\beta$, and as \cite{donoho2015higher} indicate: 
\begin{equation} \label{eqn:rho}
\rho(\beta) = \left\lbrace \begin{array}{ll} \beta - 1/2 & \text{ when } \beta \in (1/2, 3/4), \\ \left(1 - \sqrt{1 - \beta}\right)^2 & \text{ when } \beta \in (3/4, 1). \end{array}\right.
\end{equation}

If we replace the normal distribution with a $\chi_n^2$ distribution, the same results hold for tests of the form
\begin{equation} \label{eqn:HC_chi_test}
\begin{array}{ll}
H_{o, i}: & W_i \sim {\chi_n^2},  \\
H_{1, i}: & W_i \sim {\chi_n^2}(\delta),
\end{array}
\end{equation}
where $\delta$ is a non-centrality parameter and we consider $r \in (0, 1)$ such that $\delta = 2 r \log p$. 

\begin{remark} 
While Higher Criticism is typically formulated for the case of identical non-null means or parameters (all of the non-zero $\mu_i$ are identical), this constraint is not mandatory \cite{arias2011global, hall2010innovated}. Indeed, the results hold without modification for the Gaussian model with non-null means of size $\mu_i = \alpha_i \sqrt{2 \log p}$, where $\alpha_i$ is a non-negative random variable with the property that $\mathbb{P}(\alpha_i \leq \sqrt{r}) = 1$ and $\mathbb{P}(\alpha_i > \sqrt{r} - \epsilon) > 0$ for all $\epsilon > 0$ \cite{hall2010innovated}. The case of a $\chi_n^2$ distribution is similar. 
\end{remark}

As a point of interest, the test in (\ref{eqn:HC_norm_test}) can be extended to (and potentially strengthened in) the case where the $p$ tests are correlated, i.e., when the additive Gaussian noise has a non-identity covariance \cite{hall2010innovated}.

\subsubsection{Application to our Problem}

Recall that for the sum-SEPCA algorithm, we formed a vector of row-sums. That is, in the equisigned setting, taking sums across the rows of $\XX$, we obtain a vector $\yy$ where $y_i = \mu_i + \sigma z_i$, with $\mu_i = (\theta u_i) \|\vv\|_1$: this situation is exactly that of a sparse mean vector embedded in Gaussian noise. Similarly, taking sums of squares across the rows of $\XX$ (as in the $\ell_2$-SEPCA algorithm) yields scaled $\chi^2_n$ distributed random variables, of which only a few have non-zero non-centrality parameters. 

With knowledge of the noise distribution, we may compute the p-values of each row statistic: these p-values are used to form the Higher Criticism statistic (\ref{eqn:HC_stat}). As in \cite{donoho2015higher}, we may adapt the original global testing problem to a selection problem. For each p-value $p_{(i)}$, we have a value $HC_{p, i}$ of the higher criticism statistic (the value that is maximized in (\ref{eqn:HC_stat})). Rejecting each null hypothesis (that the coordinate of the corresponding row is zero) when $HC_{p, i}$ is larger than the threshold $\sqrt{2 \log \log p}$ is a variable selection procedure. We refer to the procedure based on the sum statistic as HC-sum-SEPCA and that based on the sum of squares statistic as HC-$\ell_2$-SEPCA. Importantly, we note that the form of the decision boundary $\rho$ is identical to the global testing case, and that applying Higher Criticism to our row statistics is a viable global testing procedure \cite{donoho2015higher}.

\subsection{FDR-SEPCA}

In this section, we give an summary of the algorithm for uncorrelated noise and defer the general case and details to Appendix \ref{sec:fdr_details}. We continue in the same vein as in the previous section on Higher Criticism. 

We note that in the equisigned, rank-$1$ setting, coordinate selection is equivalent to the estimation of a sparse mean vector. Let $y_i = \mu_i + \sigma z_i$, where $i \in \lbrace 1, \cdots, p\rbrace$ and the vector $\zz$ of the $z_i$ is normally distributed with mean $0$ and covariance $\mathcal{I}_p$. The mean vector $\bmu$ of the $\mu_i$ is assumed to be sparse; the goal is to estimate $\bmu$. Taking sums across the rows of $\XX$, we obtain a vector $\yy$ where $y_i = \mu_i + \sigma z_i$, with $\mu_i = (\theta u_i) \|\vv\|_1$. Hence, we are in the same setting as in the previous section. 

The following penalized least squares formulation, taken from \cite{johnstone2014adaptation}, yields an estimator for $\mu$: 
\begin{equation} \label{eqn:FDR_prob_simple}
\widehat{\bmu} = \arg \min_{\bmu} \|\yy - \bmu\|_2^2 + \sigma^2 \text{pen}\left(\|\bmu\|_0\right),
\end{equation}
where $\text{pen}(k)$ is defined as 
\begin{equation}
\text{pen}(k) = \zeta k \left(1 + \sqrt{2 \log(\nu  p / k)}\right)^2,
\end{equation}
with $\zeta > 1$; we may take $\zeta = 1 + o(1)$. The parameter $\nu$ is no smaller than $e$. We define $\|\bmu\|_0$ to be the number of non-zero coordinates of $\mu$. 

The solution to (\ref{eqn:FDR_prob_simple}) is given by hard-thresholding. Let $|y|_{(i)}$ be the $i^{th}$ order statistic of $|y_i|$, namely $|y|_{(1)} \geq \cdots \geq |y|_{(p)}$. Then if
\begin{equation}
\widehat{k} = \arg \min_{k \geq 0} \sum_{i > k} |y|_{(i)}^2 + \sigma^2 \text{pen}(k), 
\end{equation}
defining 
\begin{equation}
t_k^2 = \text{pen}(k) - \text{pen}(k - 1),
\end{equation}
the solution is to hard threshold at $t_{\widehat{k}}$. 

In this set-up, we have that 
$$t_k \approx \sqrt{\zeta} (1 + \sqrt{2 \log(\nu  p / k)}).$$ 
We provide a precise quantification of $t_k$ in Appendix \ref{sec:fdr_details}.

Hence, by computing $t_k$ and performing hard thresholding of the row sums, we can perform coordinate selection. Once again, this procedure replaces the test statistic/thresholding in Algorithm \ref{alg:gen}. 

\section{Estimation Error and Smallest Detectable Coordinate} \label{sec:risk}

As we will see, our theorems discuss the ``detectability'' of the coordinates $u_i$ of $\uu$. However, it is common in the sparse PCA literature to discuss lower bounds for the risk (estimation error) \cite{johnstone2009consistency, birnbaum2013minimax, ma2013sparse}. In what follows, we will show that these two notions are equivalent.

We define the $L^2$ estimation error for a principal component estimator as
\begin{equation} \label{eqn:loss}
L(\widehat{\uu}, \uu) = \left\|\uu - \text{sign}(\langle \uu, \widehat{\uu}\rangle) \widehat{\uu}\right\|_2^2.
\end{equation}
The quantity in (\ref{eqn:loss}) is upper bounded by $2$; this bound is attained when $\uu$ and $\widehat{\uu}$ are unit norm and mutually orthogonal. Following \cite{birnbaum2013minimax}, we want to compute a lower bound for the maximum expected loss for the $s$-sparse vectors $\uu$ (in the sense of $\ell_0$ sparsity) defined as
\begin{equation} \label{eqn:sup_loss}
\sup_{\uu \in \mathbb{S}^{p - 1} : \|\uu\|_0 \leq s} \mathbb{E} L\left(\widehat{\uu}, \uu\right),
\end{equation}
where $\mathbb{S}^{p - 1}$ denotes the unit sphere in $\mathbb{R}^p$. Let $\widehat{I}$ be some index set of coordinates selected by an algorithm of the form given in Algorithm (\ref{alg:gen}). We may take $\langle \uu, \widehat{\uu}\rangle$ to be non-negative, and decompose the loss as
\begin{equation} \label{eqn:loss_decomp}
\|\uu - \widehat{\uu}\|_2^2 = \underbrace{\|\uu_{\widehat{I}} - \widehat{\uu}\|_2^2}_{\substack{\text{Estimation Error from }\\ \text{detected coordinates}}} + \underbrace{\|\uu_{\widehat{I}^c}\|_2^2}_{\substack{\text{Error from }\\ \text{missed coordinates}}} \geq \|\uu_{\widehat{I}^c}\|_2^2.
\end{equation}

Equation (\ref{eqn:loss_decomp}) shows that the loss is lower-bounded by the squared sum of the missed coordinates. Indeed, it is a natural consequence of the result in \cite{FBG_RRN_2012} that if the sparsity $s$ grows slower than does $n$, and we have a consistent estimate of the support of $\uu$, the estimation error will asymptotically be small. Essentially, we are estimating the singular vectors of an $s \times n$ matrix instead of a $p \times n$ matrix, so that if the ratio $s / n$ has limit zero, our estimates will be consistent (see (\ref{eqn:break}) and \cite{FBG_RRN_2012}). This suggests the following strategy for lower-bounding (\ref{eqn:sup_loss}): we want to construct a non-trivial `worst-case' sparse vector. That is, we want a vector $\uu$ that has a non-trivial loss (less than $2$), is sparse (fewer than $s$ non-zero coordinates), and has maximal error from missed coordinates. To ensure a non-trivial loss, we set the first coordinate $u_1$ to be large, \textit{i.e.}, $u_1 = \sqrt{1 - r^2}$, where $r = o(1)$. To ensure sparsity, we set $u_2, \cdots, u_{m + 1}$ to be non-zero for some $m \leq s - 1$, with the subsequent coordinates of $u$ set to $0$. Then, the expected loss has the lower bound
\begin{equation}
\begin{split}
\mathbb{E} L(\uu, \widehat{\uu}) &\geq \sum_{k = 1}^p |u_k|^2 \mathbb{P}\left(\text{Not Selecting Coordinate k}\right) \\
&\geq \sum_{k = 2}^{m + 1} |u_k|^2 \mathbb{P}\left(\text{Not Selecting Coordinate k}\right),
\end{split}
\end{equation}
since $u_1$ is detected with probability approaching $1$ and $u_k$ is zero for $k > m + 1$. Now, let $u_2$ through $u_{m + 1}$ all have value ${r} / \sqrt{m}$, so that we may simplify the lower bound to
\begin{equation}
\mathbb{E} L(\uu, \widehat{\uu}) \geq r^2\mathbb{P}\left(\text{Not Selecting Coordinate k}\right).
\end{equation}
If coordinates of size $r / \sqrt{m}$ are not detected with a probability approaching $1$, $r^2$ is a lower-bound on the risk. This construction shows that specifying the sizes of coordinates that are not detected with probability approaching $1$ is equivalent to specifying a worst-case risk lower bound. {Note that the value of $r^2$ depends on the specific algorithm and estimator, and that this is not a general or universal bound. Rather, the purpose of this construction is to show the equivalence between the two perspectives (a lower bound and detectable coordinate size).}

Consequently, in what follows we focus on the smallest detectable and largest undetectable coordinates because they directly shed light on the attainable estimation error. The details of the risk calculations and extensions to approximate sparsity are deferred to Appendix \ref{sec:lq_risk}, where we summarize our findings in Theorem \ref{thm:risk}.

\section{Main Results} \label{sec:perf}

The following theorem characterizes consistent support recovery conditions. These results are the analogue of the `sparsistency' guarantees found in the LASSO and $\ell_1$-norm minimization literature \cite{ravikumar2010high}. Throughout, $\widehat{I}$ denotes the set of coordinates selected by the coordinate selection scheme. 
\begin{theorem} \label{thm:sparsistency}
For the model specified in (\ref{eqn:model}) and (\ref{eqn:u}) and the algorithms specified in Table \ref{tab:details}, assume that $p(n), n \rightarrow \infty$, $s(n) / n \rightarrow 0$, and $\log p(n) = o(n)$.
Let $\epsilon \in (0, 1)$. We have that
\renewcommand{\theenumi}{\alph{enumi}}
\begin{enumerate}[leftmargin = *]
\item{For $i \in I^{c}$,
$$\phantom{\qquad} \max_{i \in I^c} ~\mathbb{P}\left(i \in  \widehat{I}\right)~
\rightarrow 0 ,$$}
\item{For $i \in I$,
$$\min_{\substack{i \in I~:~  |\theta u_i| > \beta_{crit} (1 + \epsilon)}} 
~\mathbb{P}\left(i \in  \widehat{I}\right)~
\rightarrow 1,$$ 
$$\max_{\substack{i \in I ~:~ |\theta u_i| < \beta_{crit} (1 - \epsilon)}}
~\mathbb{P}\left(i \in  \widehat{I}\right)~ \rightarrow 0.$$ 
}
\end{enumerate}
Here
\begin{equation} \label{eqn:limits}
\beta_{crit} =
\left \lbrace \begin{array}{ll}
\sigma C_U \frac{\sqrt{\log p}}{\left| \sum_{k} v_k \right|} & \text{for sum-SEPCA,} \\ 
\sigma \sqrt{C_2} \sqrt{\frac{\log e p}{\sqrt{n}}} & \text{for $\ell_2$-SEPCA,} \\
\sigma t_{\ell_1} & \text{for $\ell_1$-SEPCA,}\\[0.15cm]
\end{array}\right.
\end{equation}
and  $t_{\ell_1}$ satisfies the relation
\begin{align*}
\left(\sqrt{\frac{2}{\pi}} + C_1 \frac{\log e p}{\sqrt{n}}\right) &= 
\frac{1}{n} \sqrt{\frac{2}{\pi}} [\sum_k \exp\left(-\left(\sqrt{n} \frac{(t_{\ell_1}) v_k}{ \sqrt{2}}\right)^2\right) + \\
&\sqrt{\pi} \sum_k \left(\sqrt{n} \frac{(t_{\ell_1})  v_k}{ \sqrt{2}}\right) \text{Erf} \left(\sqrt{n} \frac{(t_{\ell_1})  v_k}{ \sqrt{2}}\right)].
\end{align*}
\end{theorem}

We defer the proof to Appendix \ref{sec:thm1_pf}.



Theorem \ref{thm:sparsistency} identifies a phase transition in the ability of the algorithms to accurately estimate the support of $\uu$. Note that the analysis brings into sharp focus the dependence of $\beta_{crit}$ on $\vv$ for the $\ell_1$- and sum-SEPCA algorithms, but not the $\ell_2$-SEPCA algorithm. Consequently, we can expect the algorithms to perform differently depending on the structure of the underlying $\vv$. It is important to note that the sparsity $s$ of $\uu$ is not a parameter in the thresholds and results. 

It is also important to note that $\ell_2$-SEPCA and $\ell_1$-SEPCA do not rely on the equisigned character of $\vv$. However, it is clear that the sum-SEPCA algorithm explicitly depends on the equisigned assumption.

\subsection{Hamming Loss}

{It is also possible to state the above results in terms of the Hamming loss for the support of $\uu$, and prove consistency of the coordinate selection scheme by assuming that all the nonzero coordinates of $\uu$ lie above a critical threshold.
A detailed decision-theoretic analysis of variable selection under a sequence model with \emph{i.i.d.} noise and with respect to the Hamming loss has recently been carried out by  \cite{butucea2018variable}.
Recall that the Hamming loss measures the number of elements in two sets that are different, so that here the loss between the true support $I$ and the estimated support $\widehat{I}$ would be the size of the symmetric set difference of $I$ and $\widehat{I}$. Let $d_H\left(I, \widehat{I}\right)$ denote the Hamming loss and assume that whatever algorithm we are using has a threshold $\beta_{crit}$. Then, for any $\epsilon \in (0,1)$, we may write
\begin{eqnarray} \label{eqn:expected_hamming}
    \EE d_H\left(I, \widehat{I}\right) &=& \sum_{i \in I : |\theta u_i| > \beta_{crit} (1 + \epsilon)} \Prob\left(i \notin \widehat{I}\right) + \sum_{i \in I : |\theta u_i| < \beta_{crit} (1 - \epsilon)} \Prob\left(i \notin \widehat{I}\right) \nonumber\\
    && + \sum_{i \in I : (1 - \epsilon) \leq |\theta u_i|/\beta_{crit}  \leq (1 + \epsilon)} \Prob\left(i \notin \widehat{I}\right) + \sum_{i \notin I} \Prob\left(i \in \widehat{I}\right).
\end{eqnarray}
We can then restate the results on coordinate selection in terms of the Hamming loss under a more restricted setting that assumes an exact form of sparsity of the vector $\uu$.
\begin{corollary} \label{cor:hamming}
For the model specified in (\ref{eqn:model}) and (\ref{eqn:u}), and an algorithm specified in Table \ref{tab:details}, assume that the conditions of Theorem \ref{thm:sparsistency} hold. Let the support $I$ of $\uu$ be estimated by $\widehat{I}$. Moreover, assume that the algorithm has a threshold $\beta_{crit}$ (given in (\ref{eqn:limits})) such that
for a small, fixed $\epsilon_0 > 0$, the set
\begin{equation}
    I_0 := \left\{j : |\theta u_j| > \beta_{crit} (1 + \epsilon_0)\right\}
\end{equation}
equals the set $I$. 
Then the expected Hamming loss satisfies
\begin{equation}
    \EE d_H\left(I, \widehat{I}\right)/|I| \rightarrow 0.
\end{equation}
\end{corollary}
The proof of the corollary, given in Appendix \ref{ssec:cor_pf}, follows from applying Theorem \ref{thm:sparsistency},
with a more detailed enumeration of the sets and 
the inclusion probabilities,
to each term of (\ref{eqn:expected_hamming}).
}

\subsection{FDR-Based Algorithms}

We may summarize the coordinate selection properties of the FDR refinements as follows: 
\begin{theorem} \label{thm:sparsistency_fdr}
For the model specified in (\ref{eqn:model}) and (\ref{eqn:u}) and the three FDR-controlling algorithms summarized in Algorithm \ref{alg:gen_fdr}, assume that $p(n), n \rightarrow \infty$, $s(n) / n \rightarrow 0$, and 
$\log p(n) = o(n)$.
Let $\epsilon \in (0, 1)$. We have that
\renewcommand{\theenumi}{\alph{enumi}}
\begin{enumerate}[leftmargin = *]
\item{For all three algorithms and $i \in I^{c}$,
$$\phantom{\qquad} \max_{i \in I^c} 
~\mathbb{P}\left(i \in  \widehat{I}\right)~
\rightarrow 0 ,$$}
\item{For the Higher Criticism-based algorithms and $i \in I$,
$$\min_{\substack{i \in I~:~ |\theta u_i| > \beta_{crit} (1 + \epsilon)}} 
~\mathbb{P}\left(i \in  \widehat{I}\right)~\rightarrow 1,$$ 
$$\max_{\substack{i \in I ~:~|\theta u_i| < \beta_{crit} (1 - \epsilon)}}
~\mathbb{P}\left(i \in  \widehat{I}\right)~\rightarrow 0.$$ 
}
\item{For the FDR-SEPCA algorithm, uniformly over $i \in I$, \\
if $|\theta u_i| > \beta_{crit} (1 + \epsilon)$, coordinate $i$ is selected; \\
if $|\theta u_i| < \beta_{crit} (1 - \epsilon)$, coordinate $i$ is not selected\\
with probability tending to 1.}
\end{enumerate}
Here
\begin{gather} \label{eqn:limits_fdr}
\beta_{crit} =
\left \lbrace \begin{array}{ll}
\sigma \sqrt{\rho(\beta)} \frac{\sqrt{2 \log p}}{\|\vv\|_1} & \text{for HC-sum-SEPCA,} \\ 
\sigma \rho(\beta) \frac{2 \log p}{\sqrt{n}} & \text{for HC-$\ell_2$-SEPCA,} \\
\sigma \left(1 - o(1)\right) \sqrt{\zeta} \frac{1 + \sqrt{2 \log (\nu  p / \widehat{k})}}{\|\vv\|_1} & \text{for FDR-SEPCA,}\\
\end{array}\right.
\end{gather}
where $\zeta > 1$, $\nu > e$, and the FDR-SEPCA algorithm detects $\widehat{k}$ coordinates. 
\end{theorem}
We defer the proof to Appendix \ref{sec:thm2_pf}.

Once again, we see that the structure of the underlying $\vv$ plays a role in the performance of the sum-based algorithms, but not for the $\ell_2$-based HC-$\ell_2$-SEPCA algorithm. Unlike in the FWER-controlling cases, the sparsity of $\uu$ plays a (small) role here, via the constant $\rho\left(\beta\right)$ for the Higher Criticism-based methods and via $\widehat{k}$ for FDR-SEPCA. Moreover, $\ell_2$-HC-SEPCA, like $\ell_2$-SEPCA, does not make use of the equisigned nature of $\vv$. 


\subsection{Higher Ranks}

{In this work, we restrict our focus to the rank-$1$, equisigned setting. A natural question is are our results extensible to the higher rank setting?}

{The first point is concerned with the right singular vectors. To preserve orthogonality, we would need equisigned right-singular vectors $\vv_i$ with disjoint supports. The second point is concerned with the left singular vectors. Our algorithms are based on thresholding row-statistics: it is possible that the union of supports of several sparse vectors is a relatively large set. The FWER-controlling algorithms (by design) are not sensitive to the increased supports, but the FDR-controlling algorithms are sensitive to this. Indeed, the decision boundaries for the FDR algorithms explicitly depend on the sparsity levels. Third point, once again, is concerned with the left singular vectors. It is possible that a sum-based statistic suffers from cancellations that decrease the size of the row-statistic. For example, in a rank-$2$ setting, if $\uu_1 = \frac{1}{\sqrt{2}} \begin{bmatrix} 1 & 1 & 0 & \cdots & 0\end{bmatrix}^T$ and $\uu_2 = \frac{1}{\sqrt{2}} \begin{bmatrix} 1 & -1 & 0 & \cdots & 0\end{bmatrix}^T$ and $\left\|\vv_1\right\|_1$ and $\left\|\vv_2\right\|_1$ have similar values and are both non-negative, the row-sum of the second row will be small. Note, however, that the $\ell_2$-norm based methods do not suffer from this issue.}

\section{Simulations} \label{sec:sim}

To illustrate the relative powers of the six algorithms, we compute the theoretical limits on the sizes of detectable coordinates as a function of $n$. We use a unit-norm, equisigned $\vv$ such that
\begin{equation} \label{eqn:v_good}
v_k \propto \exp \left(-5 \frac{k}{n}\right) \left|\sin\left( 4 \frac{k}{n} \right)\right| \text{ for } 1 \leq k \leq n.
\end{equation}

\begin{figure}
\begin{minipage}{\linewidth}
\centering
\includegraphics[width = 0.85\textwidth]{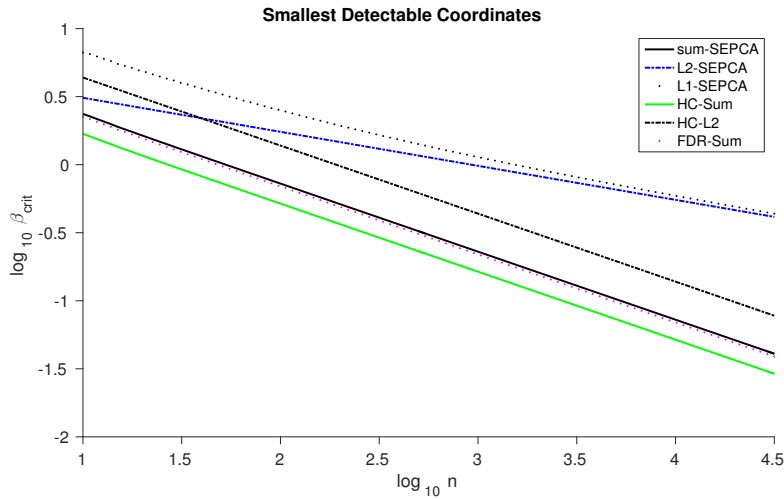}
\captionof{figure}{This plot shows $\beta_{crit}$ for all six algorithms for the $\vv$ described in (\ref{eqn:v_good}).}
\label{fig:size}
\end{minipage}
\end{figure}

This choice of $\vv$ has a `rise and fall' sort of behavior, and is motivated by physical signals, e.g., chemical reactions or nerve signals in the brain. The value of $\beta_{crit}$ is shown in Figure \ref{fig:size}; for this choice of $\vv$, it is clear that the sum-SEPCA dramaticaly outperforms the other SEPCA variants in terms of size of the smallest detectable component. The FDR-SEPCA algorithm has similar performance to sum-SEPCA, and the HC-sum-SEPCA algorithm has the strongest performance. 

\begin{figure}
\begin{minipage}{\linewidth}
\centering
\includegraphics[width = 0.95\textwidth]{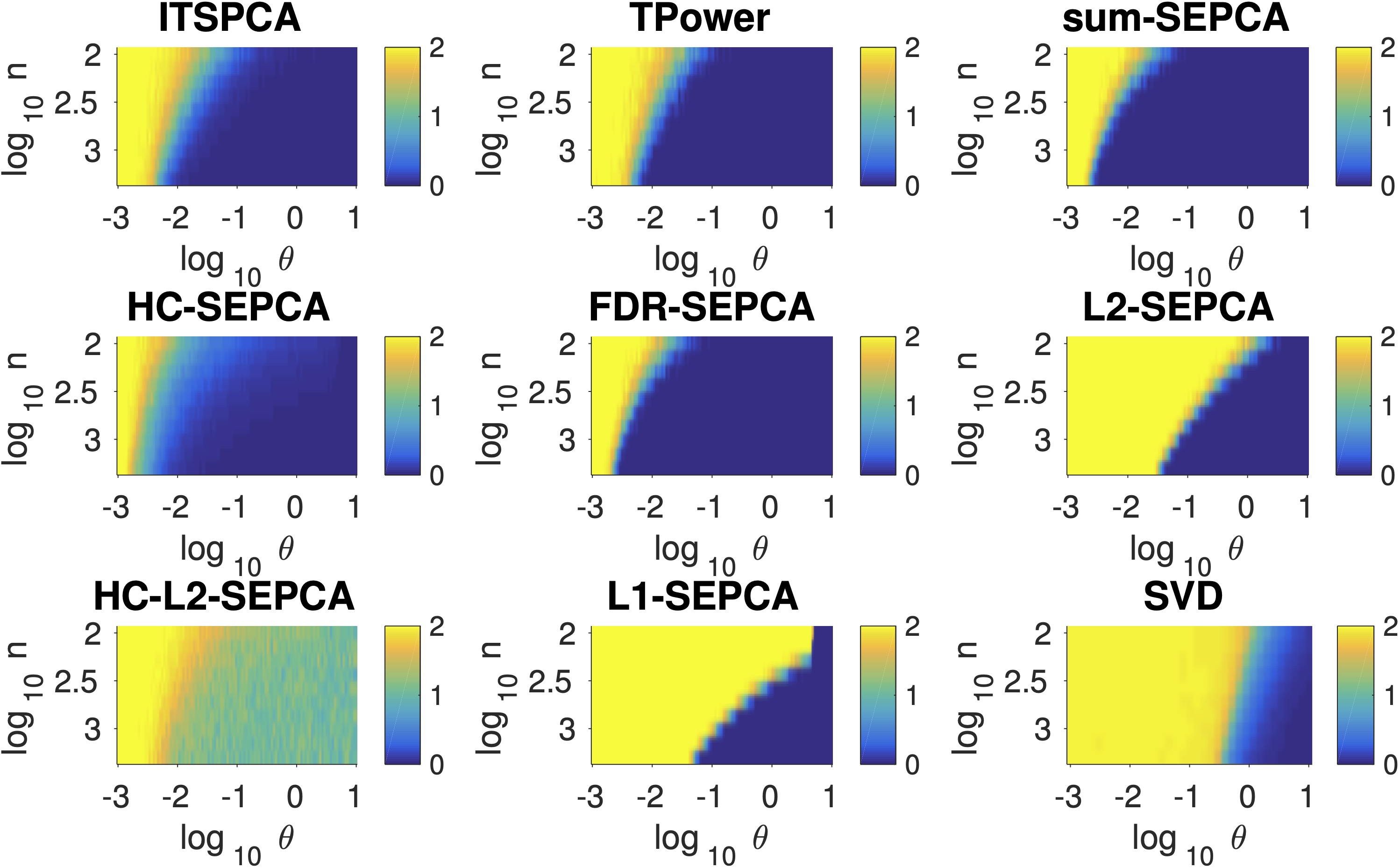}
\captionof{figure}{The plots show the empirical estimation error for all six algorithms for the $\uu$ with one non-zero coordinate and the $\vv$ described in  (\ref{eqn:v_good}). We include results from TPower, ITSPCA and the SVD for comparison.}
\label{fig:risk}
\end{minipage}
\end{figure}

In Figure \ref{fig:risk}, we plot the estimation error as a function of $n$ and $\theta$ for all six algorithms. We also include results for the SVD and competing algorithms TPower \cite{yuan2013truncated} and ITSPCA \cite{ma2013sparse}. In the simulations, we fix $p = 1000$ and vary $n$, since the dependence in $p$ in the thresholds is logarithmic, whereas that in $n$ is not. The left singular vector $\uu$ is chosen to be the vector with $1$ in the first coordinate and $0$ elsewhere. We fix the noise variance $\sigma^2$ at $1$, so that $\theta^2$ is the eigen-SNR. The results should be interpreted as follows. For the particular $\vv$ chosen here, we expect HC-sum-SEPCA to have the lowest detectable limit, and $\ell_1$-SEPCA to have the largest. This behavior is confirmed. Moreover, the sum-based algorithms offer a slight strengthening of both ITSPCA and TPower. Importantly, note that the sum-based algorithms explicitly take advantage of the equisigned nature of $\vv$: that is, algorithms that explicitly use the equisigned property outperform algorithms that do not (the $\ell_2$ and $\ell_1$ algorithms, as well as ITSPCA and TPower). 

\noindent

\begin{figure}
\begin{minipage}{\linewidth}
\centering
\includegraphics[width = 0.95\textwidth]{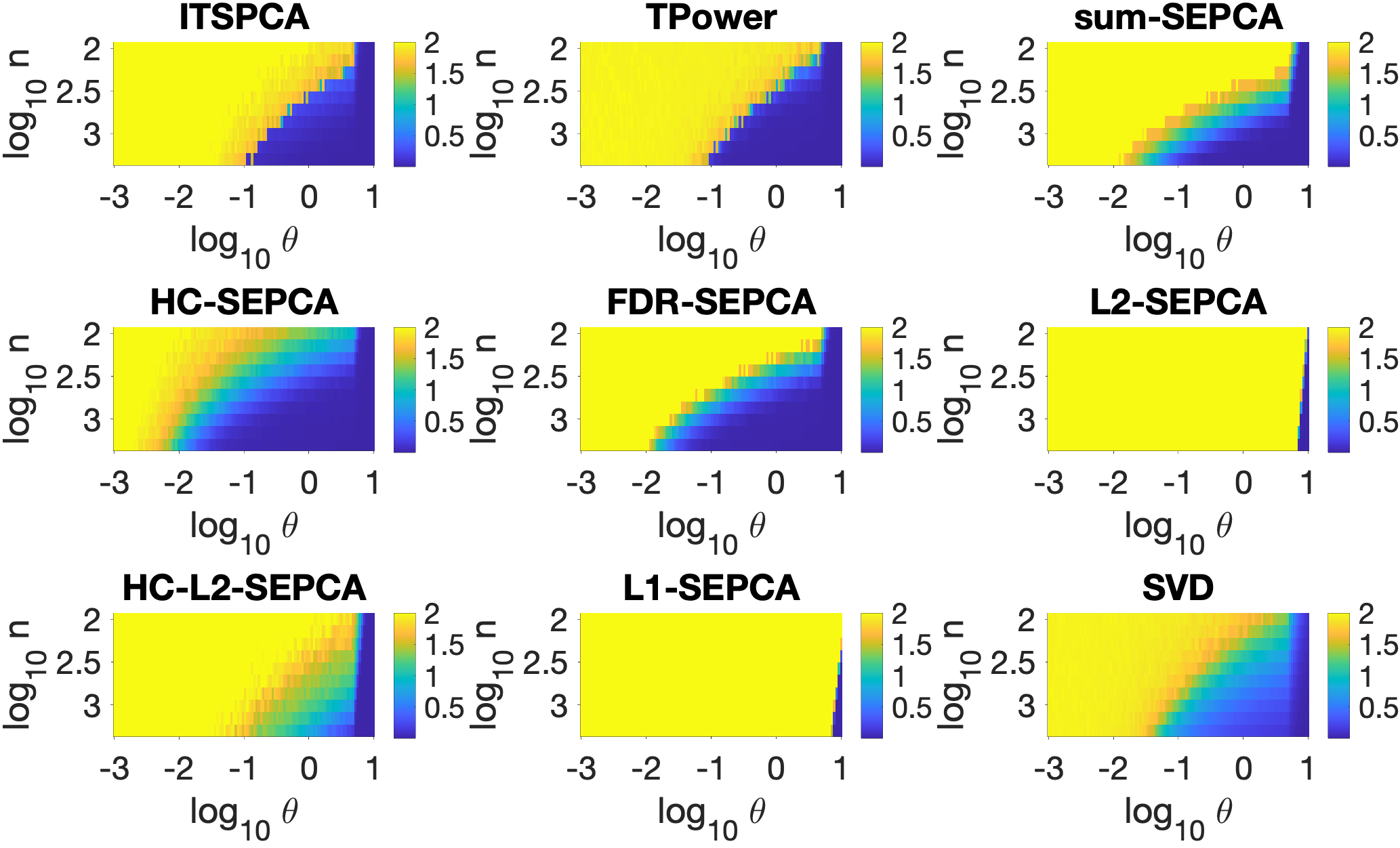}
\captionof{figure}{The plots show the empirical estimation error for all six algorithms for the $\uu$ with $\sqrt{p}$ non-zero coordinates and the $\vv$ described in  (\ref{eqn:v_good}). We include results from TPower, ITSPCA and the SVD for comparison.}
\label{fig:risk_dense}
\end{minipage}
\end{figure}

{We repeat our simulations for a $\uu \in \RR^p$ with $\sqrt{p}$ non-zero coordinates (of equal size) and the same $\vv$, as seen in Figure \ref{fig:risk_dense}. We find similar conclusions as in Figure \ref{fig:risk}, where the sum-based algorithms offer a strengthening over ITSPCA and TPower; the $\ell_2$-norm based algorithms do not perform well. Note that a sparsity of $\sqrt{p}$ is at the limit/valid edge for the higher criticism-based methods, but that these methods still perform well.}

\subsection{Comments on the FDR-controlling procedures }
%

The Higher Criticism for the $\chi_n^2$-variates `pushes back' the phase transition between detecting nothing and something to a lower value of $\theta$ relative to the $\ell_2$-SEPCA algorithm, but is still less powerful than any of the sum-based algorithms. Moreover, even above the phase transition, the $\ell_2$-SEPCA algorithm may be preferable, as the error is increased by unacceptably many false positives. 

The Higher Criticism procedure for the sum statistic has the lowest phase transition point and hence the highest power. Its transition is more gradual than the 
penalized FDR thresholding procedure and sum-SEPCA, which have roughly the same performance in this simulation. 

\subsection{An example where $\ell_2$-based algorithms outperform sum-based algorithms}

\begin{figure}[htb] 
\begin{minipage}[b]{\linewidth}
\centering
\includegraphics[width = 0.85\textwidth]{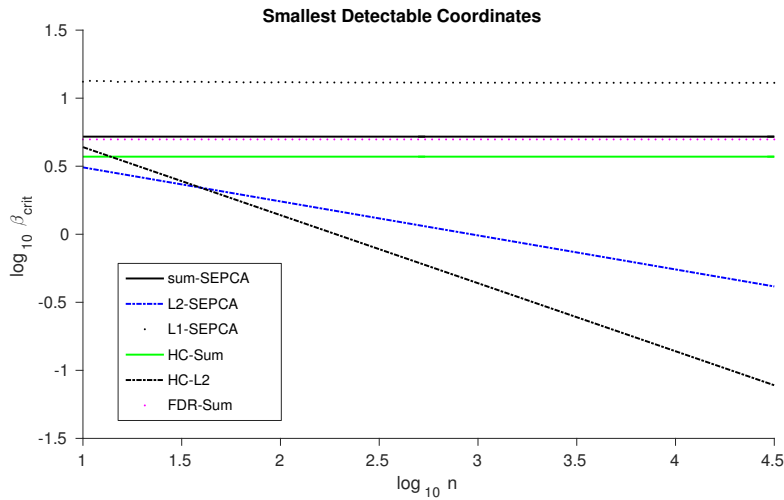}
\captionof{figure}{This plot shows $\beta_{crit}$ for all six algorithms for the $\vv$ described in (\ref{eqn:v_bad}).}
\label{fig:size_bad}
\end{minipage}
\end{figure}

\begin{figure}[htb] 
\begin{minipage}[b]{\linewidth}
\centering
\includegraphics[width = 0.95\textwidth]{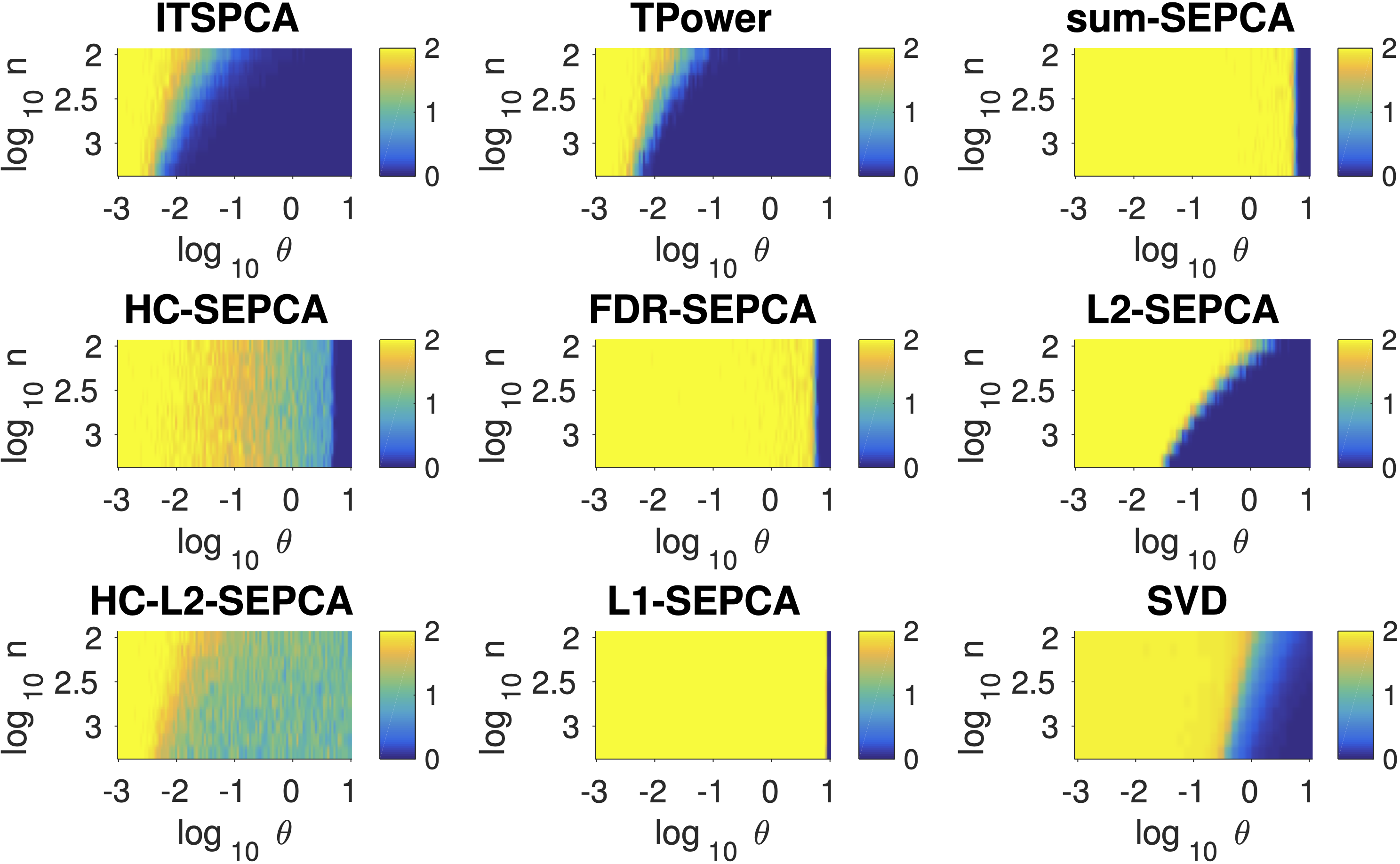}
\captionof{figure}{The plots show the empirical estimation error for all six algorithms for the $\uu$ and $\vv$ described in  (\ref{eqn:v_bad}). We include results from TPower, ITSPCA and the SVD for comparison.}
\label{fig:risk_bad}
\end{minipage}
\end{figure}

Sum-SEPCA has a $\beta_{crit}$ that depends on $\vv$. Looking at the form in (\ref{eqn:limits}), if $\|\vv\|_1$ is smaller than $n^{1/4}$, we would expect $\ell_2$-SEPCA to detect a smaller coordinate size. Vectors with smaller coordinates have a smaller $\ell_1$-norm, i.e., one that is closer to their $\ell_2$-norm. Hence, if we choose 
\begin{equation} \label{eqn:v_bad}
v_k \propto \frac{1}{k^2} \text{ for } 1 \leq k \leq n,
\end{equation}
we expect sum-SEPCA to have worse performance relative to $\ell_2$-SEPCA. Figures \ref{fig:size_bad} and \ref{fig:risk_bad} confirm this expectation. The FDR refinements perform poorly. It should be noted, however, that TPower and ITSPCA retain their performance. This choice of $\vv$ effectively corresponds to a very small value of $n$: the majority of coordinates are tiny in size and buried beneath noise regardless of the value of $\theta$. If we `corrected' the scenario and used a smaller $n$ and a subset of $\vv$, we would be in a situation closer to that given in (\ref{eqn:v_good}).

\subsection{A video data example}

We conclude our sequence of examples with a real data study. This example is motivated by the problem of foreground-background separation in videos. Consider a grayscale video of stars twinkling against a black background \cite{TSvid}. Our goal is to estimate the locations of the stars: by reshaping the video, we may treat each frame as a vector and hence treat the video as a sparse matrix. Only a few locations have a star and are hence non-zero. The scale of the video pixels is between $0$ and $255$. We examine the top-left $72\times 64$ pixels for $89$ frames, as shown in Figure \ref{fig:TS_vid_im}. In Figure \ref{fig:TS_vid_sing}, we plot the singular values of the video matrix. The first singular value stands out strongly against the rest, and at most two more singular values are well-separated from the bulk. This structure suggests that our rank-$1$ based approach is well suited to this problem. 

\begin{table}[ht]
\begin{minipage}{\linewidth}
\begin{tabular}{cc}
\begin{subfigure}{0.45\textwidth}\centering\includegraphics[width=\columnwidth, trim = {0 0 0 1cm}, clip]{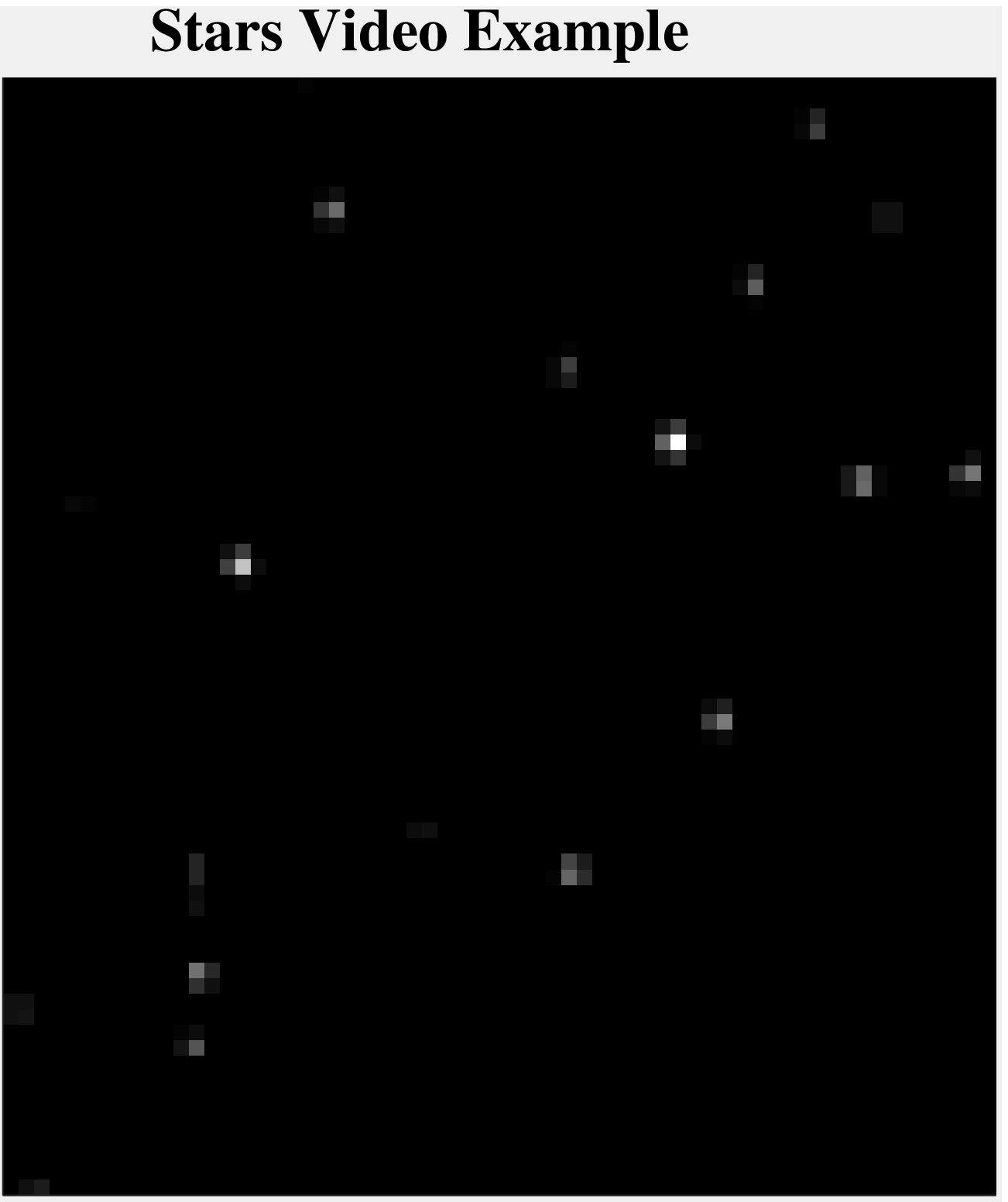}\caption{The image shows the mean intensity of pixels from the top-left $72\times 64$ pixels for $89$ frames. White indicates the presence of a star.}\label{fig:TS_vid_im}\end{subfigure} & 
\begin{subfigure}{0.45\textwidth}\centering\includegraphics[width=\columnwidth]{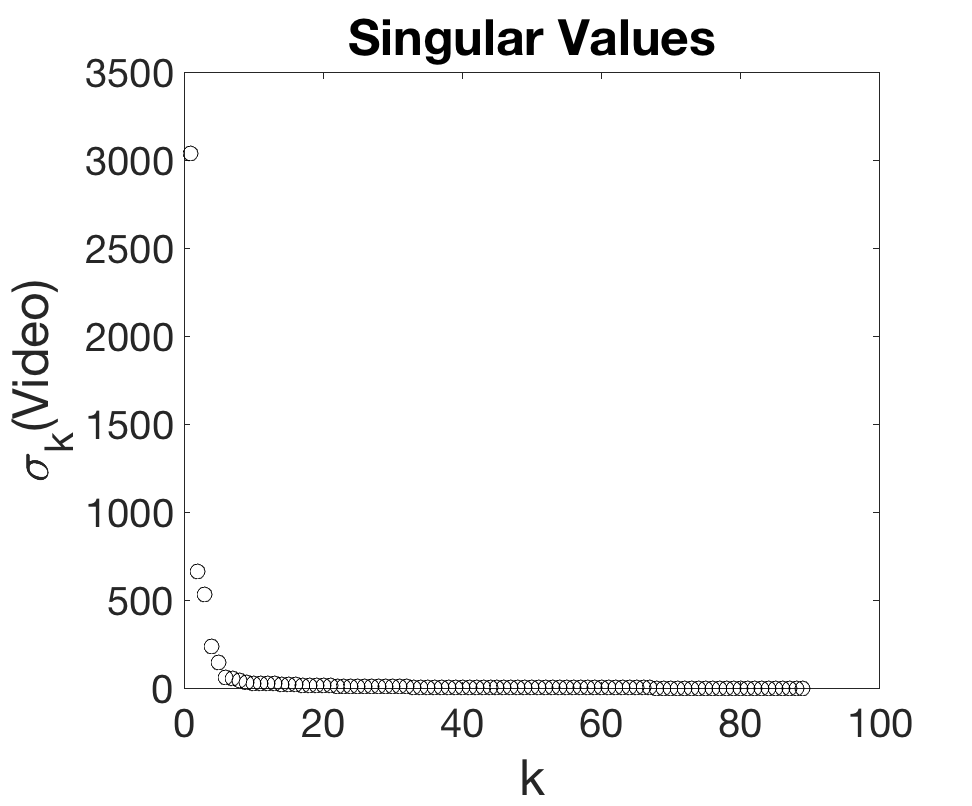}\caption{The plots shows the singular values of the video data. The spacing suggest a low-rank-plus-noise structure. }\label{fig:TS_vid_sing}\end{subfigure}\   
\end{tabular}
\caption{Video Example Figures}
\label{tab:vidfigs}
\end{minipage}
\end{table}

\begin{figure}
\begin{minipage}{\linewidth}
\centering
\includegraphics[width = 0.95\textwidth]{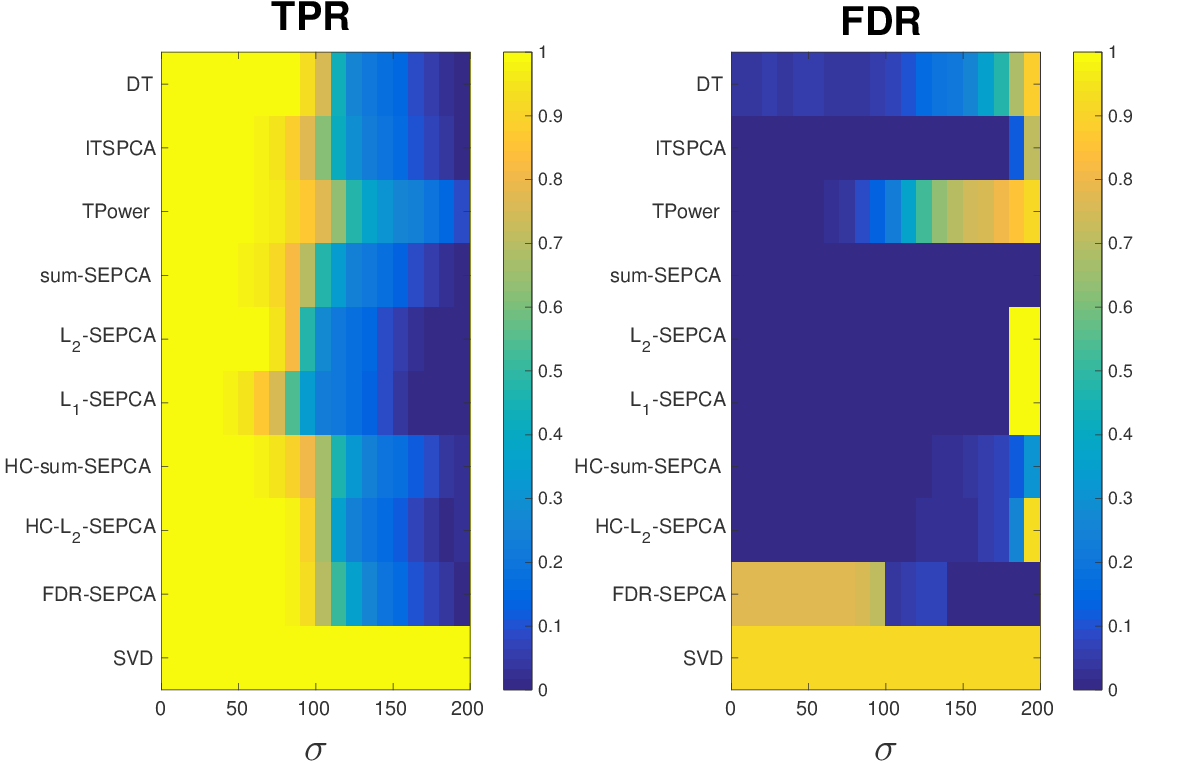}
\captionof{figure}{The left plot shows the True Positive Rate of the various algorithms as a function of the noise level $\sigma$. The right plot shows the False Discovery Rates. }
\label{fig:TS_vid_perf}
\end{minipage}
\end{figure}

\begin{figure}[htb]
\begin{minipage}{\linewidth}
\begin{subfigure}{0.49\textwidth}\centering\includegraphics[width=\columnwidth]{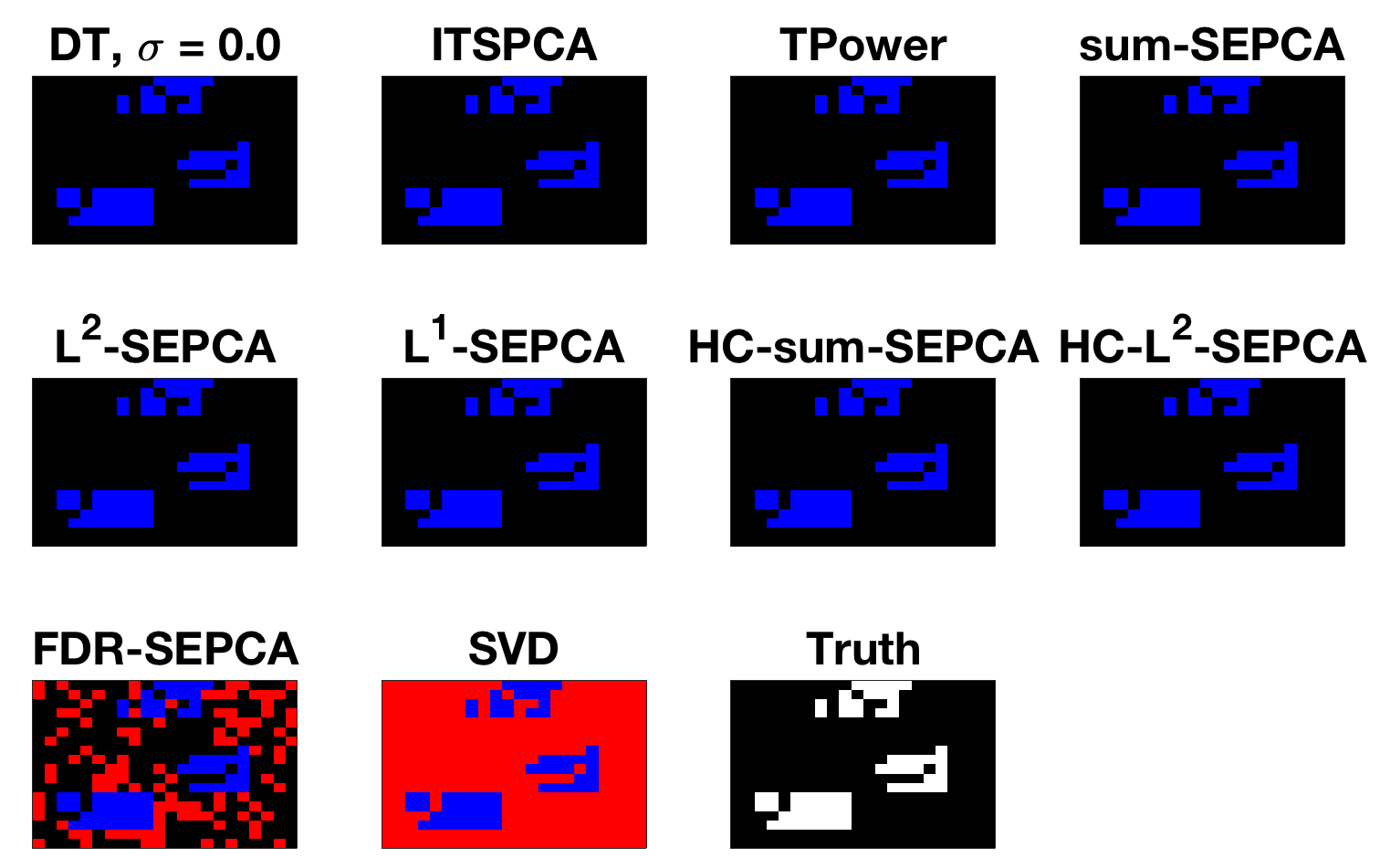}\caption{$\sigma=0$}\label{fig:TSVid_s0}\end{subfigure} 
\begin{subfigure}{0.49\textwidth}\centering\includegraphics[width=\columnwidth]{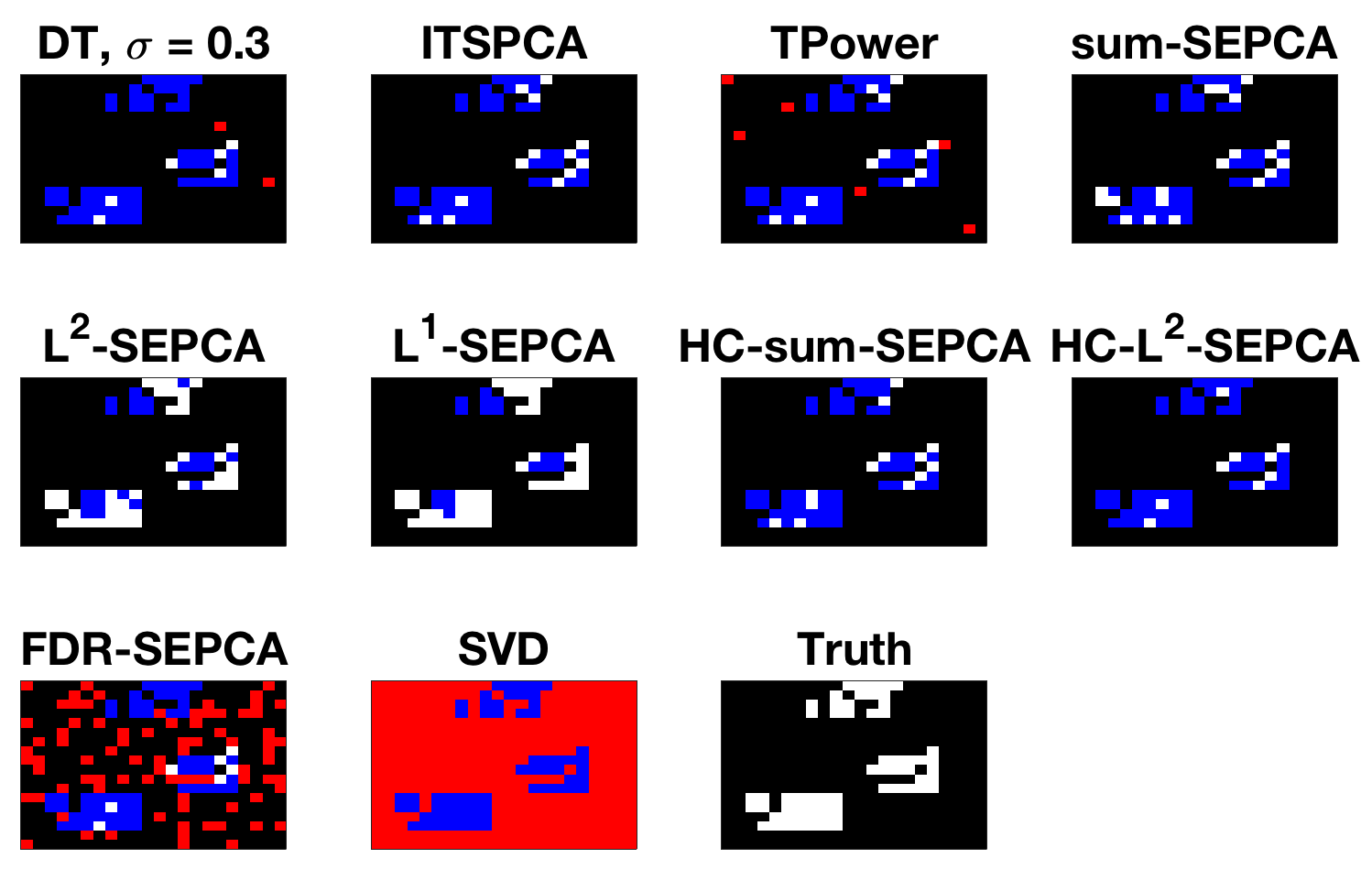}\caption{$\sigma=0.3$}\label{fig:TSVid_s03}\end{subfigure}\\
\begin{subfigure}{0.49\textwidth}\centering\includegraphics[width=\columnwidth]{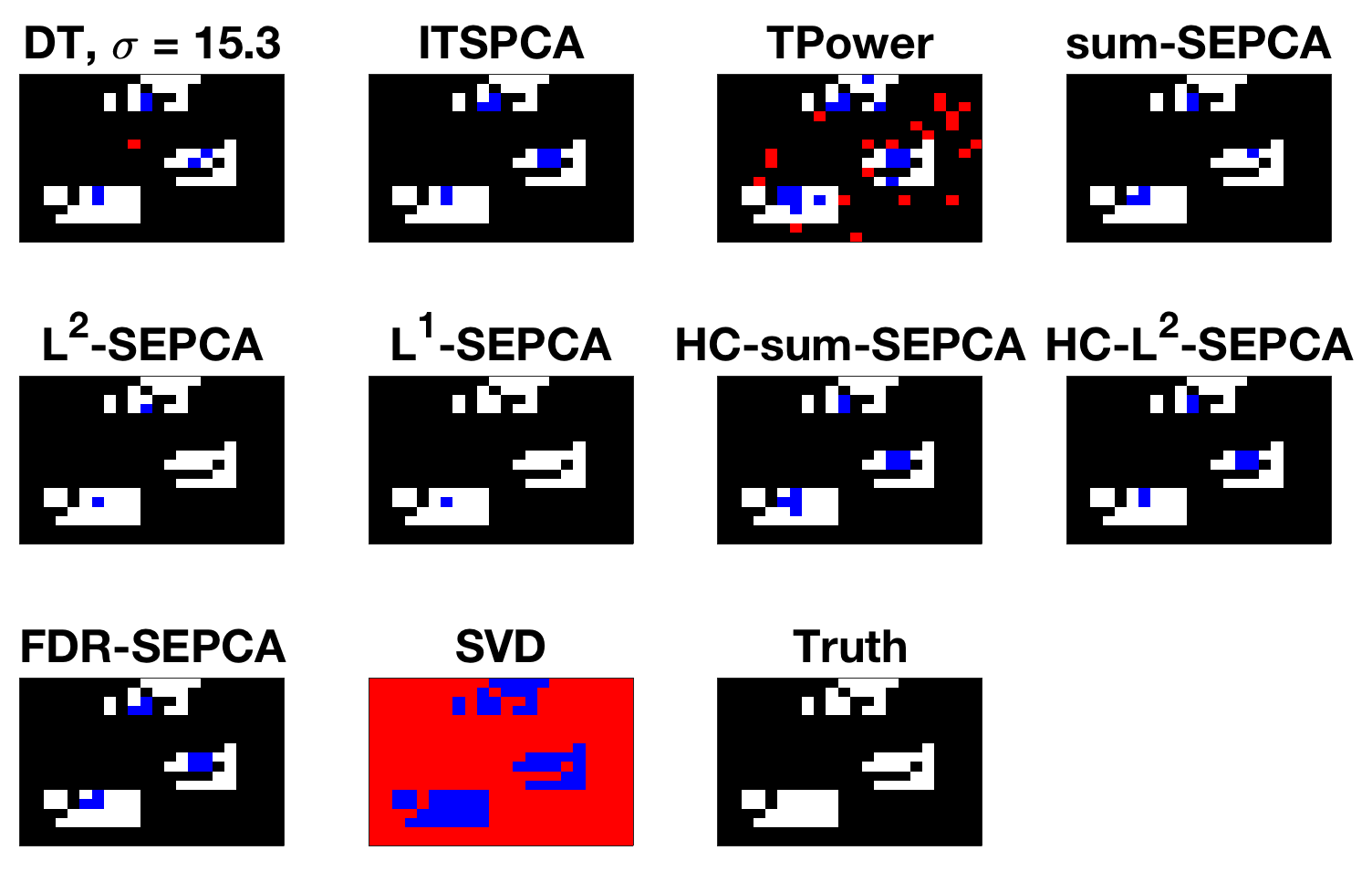}\caption{$\sigma=15$}\label{fig:TSVid_s15}\end{subfigure} 
\begin{subfigure}{0.49\textwidth}\centering\includegraphics[width=\columnwidth]{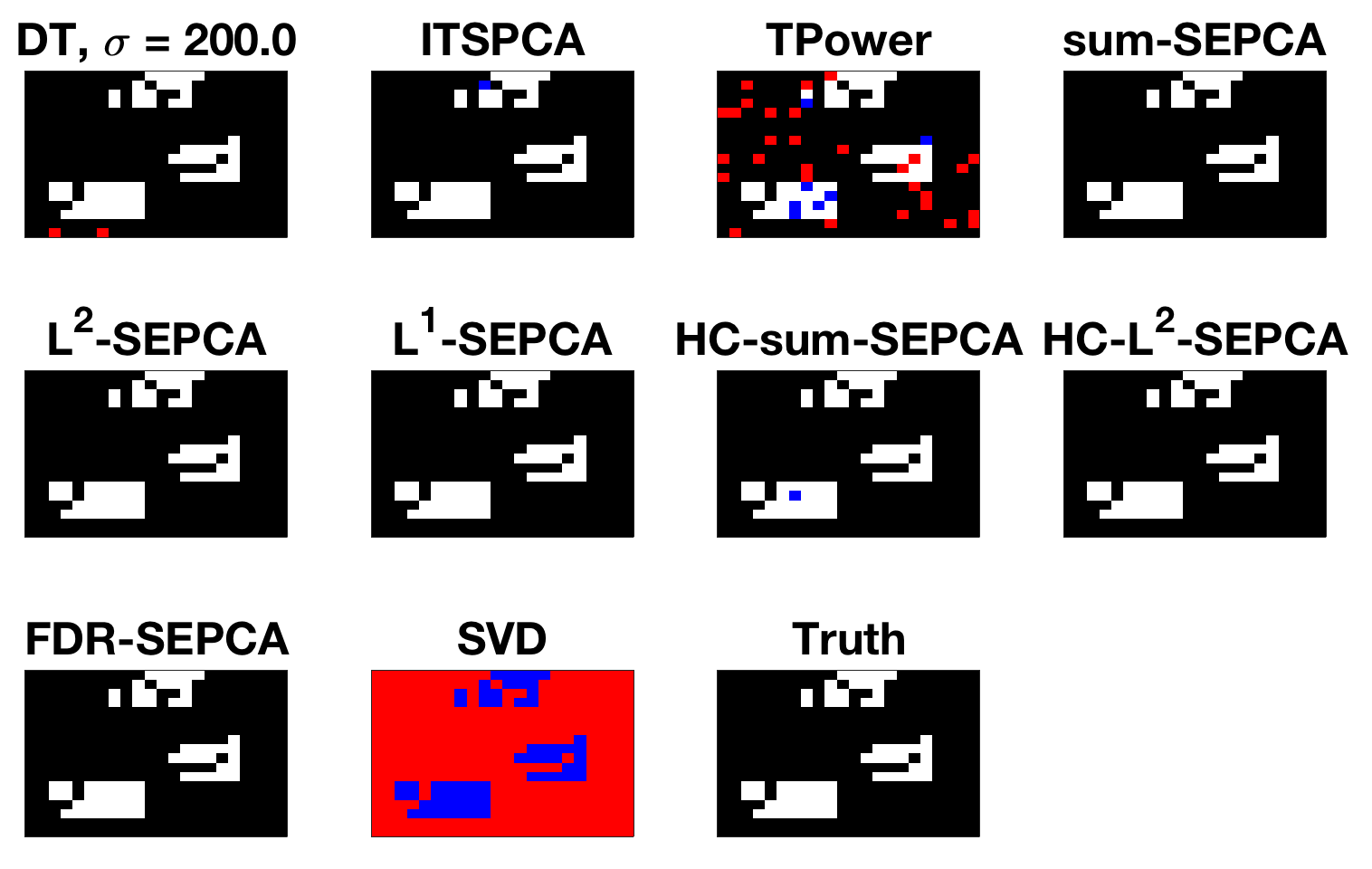}\caption{$\sigma=200$}\label{fig:TSVid_s200}\end{subfigure}\
\caption{A zoomed-in view of the three top-right stars in the video example. White indicates a false negative (missed star), Red a false positive (a guessed pixel where there was nothing), and Blue a true positive (correctly identified pixel).}
\label{tab:vidfigs_detect}
\end{minipage}
\end{figure}

We add Gaussian noise of variance $\sigma^2$ and study the True Positive Rates (TPR) and False Discovery Rates (FDR) across all algorithms and across different values of $\sigma$. In Figure (\ref{fig:TS_vid_perf}), we show the results of our simulations. In terms of the TPR, everything other than the SVD has a similar performance, while the test-statistic SEPCA-based algorithms enjoy the best performance in terms of the FDR. In Figure \ref{tab:vidfigs_detect} we zoom in on the top-right three stars and show how the algorithms perform as noise increases. Here, we see that the behavior alluded to in the TPR/FDR results actually occurs in the video. 

\section{A geometric view: which algorithm to use?} \label{sec:comp}

We have stated detectability results for each algorithm in Section \ref{sec:perf} and provided a numerical verification and comparison in Section \ref{sec:sim}. In this section, we wish to analytically compare the algorithms. In particular, we have seen that the right singular vector $\vv$ plays a critical role in the detectability and estimability of $\uu$, and we will characterize this behavior carefully. 

In this section, will use the following notational convenience: we absorb $(\theta u_i)$ into $\vv \in \mathbb{R}^n$, and write the detectability of coordinates in terms of $\vv$. That is, if $\vv^T$ is a row of $\XX$, we specify when that row is selected. Moreover, we take $\sigma = 1$ for simplicity. 

There are two `classes' of detectability: in terms of $\|\vv\|_1$ and in terms of $\|\vv\|_2$. The sum-, HC-sum, and FDR-SEPCA algorithms select a coordinate if $\left|\sum_k v_k\right| = \|\vv\|_1$ is large enough for a $\vv$ in the orthant with all non-negative or all non-positive coordinates. Geometrically, the vector $\vv$ is selected if it is `outside' a hyperplane with a normal vector proportional to the vector of all $1$s. The $\ell_1$-SEPCA algorithm is similar, as it selects a coordinate when $\|\vv\|_1$ is large enough, or if $\vv$ lies outside an $\ell_1$-ball of some radius. The connection between the previous three algorithms and $\ell_1$-SEPCA comes from noting that the faces of an $\ell_1$-ball are sections of hyperplanes with normal vectors proportional to a vector of $\pm 1$s. Finally, the $\ell_2$- and HC-$\ell_2$-SEPCA algorithms select a coordinate when $\|\vv\|_2$ is large enough. I.e., when $\vv$ lies outside some $\ell_2$-ball. 

Our goal in this section is to derive comparisons between the six algorithms. Specifically, for a given vector $\vv$, which algorithm will have the greatest detection ability (we are, for the moment, only concerned with maximizing power)? Note that when $\vv$ has a large norm, it does not matter which algorithm is used. Questions only arise when $\|\vv\|_1$ or $\|\vv\|_2$ are relatively small and are close to the thresholds. 

\subsection{Intersection of a hyperplane and a hypersphere} \label{ssec:intersect}

We may think of the $\ell_1$ ball as a hyperplane when restricted to a single orthant. If a hypersphere of radius $r$ intersects a hyperplane with a normal vector proportional to the vector of all $\pm 1$s and minimum distance to the origin of $r - h$, a hyperspherical cap of height $h$ is formed: see Figure \ref{fig:cap} for a simple illustration. Geometrically, a right triangle is formed, with hypotenuse $r$ and leg $r - h$. Hence, the angle between the center of the cap and the edge is:
\begin{equation}
\theta_{lim} = \cos^{-1} \frac{r - h}{r}.
\end{equation}
It is sufficient to guarantee that 
$$0 \leq \frac{r - h}{r} \leq 1$$ 
for the hyperspherical cap to exist. Moreover, a vector $\vv$ has a direction contained inside the cap when the angle between $\vv$ and the vector of $\pm 1$ in the orthant containing $\vv$ is smaller than $\theta_{lim}$. In other words, defining the angle for a vector $\vv$ as
\begin{equation}
\theta(\vv) = \cos^{-1} \frac{\|\vv\|_1}{\|\vv\|_2 \sqrt{n}},
\end{equation}
we need $\theta(\vv) \leq \theta_{lim}$. 

\vspace{1mm}
\noindent
\begin{minipage}[b]{\linewidth}
\centering
\includegraphics[width =0.5\textwidth]{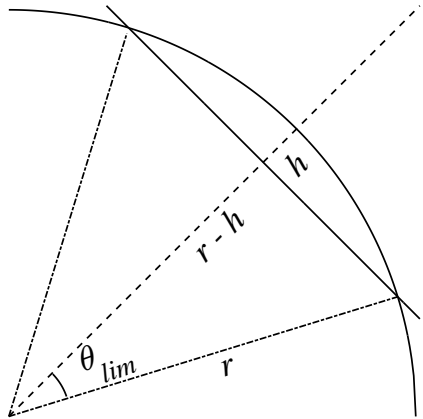}
\captionof{figure}{A spherical cap in $\mathbb{R}^2$}
\label{fig:cap}
\end{minipage}

\subsection{Comparison: $\ell_2$-based versus sum-based statistics}

We begin with a summary of the performance of each individual algorithm in Table \ref{tab:comp}. We first compare $\ell_2$-SEPCA and then compare HC-$\ell_2$-SEPCA with sum-, HC-sum, FDR-SEPCA in Tables \ref{tab:l2} and \ref{tab:HCl2}. In our comparisons, we consider when the hyperspherical cap exists and give the angle of the cap. These are routine calculations, so we omit the details. We also omit $\ell_1$-SEPCA from our comparisons, as we lack a closed-form expression for $t_{\ell_1}$. 

Note that the existence of this cap is a proxy for the equivalent statement that there exist vectors for which the sum-based algorithms are more powerful than the $\ell_2$-based algorithm. While this existence is not the same as attributing uniformly greater power to the sum-based algorithms relative to the $\ell_2$-based algorithm, the cap \emph{not} existing is equivalent to the $\ell_2$-based algorithm having uniformly greater power. 

Essentially, we observe that for $n$ and $p$ sufficiently large, the cap will exist. Moreover, for $\vv$ that is sufficiently dense ($\|\vv\|_1$ is sufficiently large), $\theta_{lim}$ will lie inside the cap. Hence, in these situations we would prefer a sum-based algorithm over an $\ell_2$-based algorithm. 

\noindent
\begin{minipage}{\linewidth} 
\captionof{table}{A summary of the six algorithms.}
\label{tab:comp}
\centering
\begin{tabular}{|l|l|l|}
\hline
Algorithm & Threshold & Geometric Quantity \\
\hline
sum & $\|\vv\|_1 \geq C_U \sqrt{\log p}$ & $r - h = C_U \sqrt{\frac{\log p}{n}}$ \\
\hline
HC-sum & $\|\vv\|_1 \geq \sqrt{2 \rho(\beta) \log p}$ & $r - h = \sqrt{2 \rho(\beta) \frac{\log p}{n}}$ \\
\hline
FDR & $\begin{aligned} &\|\vv\|_1 \geq  \left(1 - o(1)\right) \sqrt{\zeta}\\ &\left(1 + \sqrt{2 \log \left(\nu  p / \widehat{k}\right)}\right)\end{aligned}$ & $\begin{aligned}r &- h = \left(1 - o(1)\right) \sqrt{\zeta} \\ &\frac{1 + \sqrt{2 \log \left(\nu  p / \widehat{k}\right)}}{\sqrt{n}}\end{aligned}$ \\
\hline
$\ell_1$ & $\|\vv\|_1 \geq t_{\ell_1}$ & $r - h = t_{\ell_1}$ \\
\hline
$\ell_2$ & $\|\vv\|_2 \geq \sqrt{C_2} \sqrt{\frac{\log e p}{\sqrt{n}}}$ & $r = \sqrt{C_2} \sqrt{\frac{\log e p}{\sqrt{n}}}$ \\
\hline
HC-$\ell_2$ & $\|\vv\|_2 \geq 2 \rho(\beta) \frac{2 \log p}{\sqrt{n}}$ & $r = 2 \rho(\beta) \frac{2 \log p}{\sqrt{n}}$\\
\hline
\end{tabular}
\end{minipage}

\noindent
\begin{minipage}{\linewidth} 
\captionof{table}{The relative performance of $\ell_2$-SEPCA.}
\label{tab:l2}
\centering
\begin{tabular}{|l|l|l|}
\hline
Algorithm & $\cos \theta_{lim}$ & Cap exists if  \\
\hline
sum & $C_U/\sqrt{C_2} \sqrt{\frac{\log p}{(1 + \log p) \sqrt{n}}}$ & $n \geq C_U^4 / C_2^2, p \geq 1$  \\
\hline 
HC-sum & $\sqrt{2 \rho(\beta) / C_2} \sqrt{\frac{\log p}{(1 + \log p) \sqrt{n}}}$ & $n \geq 1, p \geq 1$ \\
\hline
FDR-sum & $\frac{1 + \sqrt{2 \log \nu  p/\widehat{k}}}{\sqrt{C_2 \sqrt{n} \log \nu p }}$ & $p \geq 11, n \geq 1$  \\
\hline
\end{tabular}
\end{minipage}

\noindent
\begin{minipage}{\linewidth} 
\captionof{table}{The relative performance of HC-$\ell_2$-SEPCA.}
\label{tab:HCl2}
\centering
\begin{tabular}{|l|l|l|}
\hline
Algorithm & $\cos \theta_{lim}$ & Cap exists if \\
\hline
sum & $\frac{C_U}{2 \rho(\beta) \sqrt{\log p}} $ & $p \geq \exp\left(\frac{C_U^2}{4 \rho(\beta)^2}\right)$ \\
\hline 
HC-sum & $\left[2 \rho(\beta) \log p\right]^{-1}$ & $p \geq \exp\left(\frac{1}{2 \rho(\beta)}\right)$  \\
\hline
FDR-sum & $\frac{1 + \sqrt{2 \log \nu p/\widehat{k}}}{2 \rho(\beta) \log p}$ & $\begin{aligned} &\log p \geq \\&\frac{1}{4 \rho(\beta)^2} \bigl(1 + 2\rho(\beta)\\ &+ \sqrt{8 \rho(\beta)^2 + 4 \rho(\beta) + 1}\bigr)\end{aligned}$ \\
\hline
\end{tabular}
\end{minipage}
\vspace{1mm}

\subsection{HC-$\ell_2$-SEPCA versus $\ell_2$-SEPCA}

Now, we consider when HC-$\ell_2$-SEPCA is more powerful than $\ell_2$-SEPCA. The ratio of the radii is given by 
\begin{equation} \label{eqn:HCl2_base}
\frac{2 \rho(\beta)}{\sqrt{C_2} n^{1/4}} \frac{\log p}{\sqrt{1 + \log p}}.
\end{equation}
If this ratio is smaller than $1$, HC-$\ell_2$-SEPCA is more powerful than $\ell_2$-SEPCA. Note that the quantity 
$$\frac{2 \sqrt{2}}{\sqrt{C_2}} \sqrt{\frac{\log p}{\sqrt{n}}}$$
is an upper bound for (\ref{eqn:HCl2_base}), so that if  
$$\frac{\log p}{\sqrt{n}} < \frac{e}{ 4\sqrt{2}},$$ 
the original ratio is smaller than $1$ and HC-$\ell_2$-SEPCA is preferable to $\ell_2$-SEPCA. 

\subsection{Comparing the sum-based algorithms}

Finally, we compare sum-, HC-sum-, and FDR-SEPCA. First, the ratio of the thresholds for HC-sum- and sum-SEPCA is 
\begin{equation}
\frac{\sqrt{2 \rho(\beta)}}{C_U}.
\end{equation}
Noting that $\rho(\beta) \leq 1$ and that $C_U \geq \sqrt{2} + 1/3\sqrt{2}$, it is clear that this ratio is always smaller than $1$ so that HC-sum-SEPCA is a strict improvement on sum-SEPCA.

Next, we compute the ratio of the thresholds for FDR- and sum-SEPCA: 
\begin{equation}
\frac{1 + \sqrt{2 \log \nu p/\widehat{k}}}{C_U \sqrt{\log p}}.
\end{equation}
Using the lower bound on $C_U$, we find that if $\widehat{k} \geq 11$ (and $p \geq \widehat{k}$, naturally), FDR-SEPCA is always more powerful than sum-SEPCA. For smaller values of $\widehat{k}$, for sufficiently large values of $p$, the ratio will be smaller than $1$. 

Lastly, we compare FDR-SEPCA to HC-sum-SEPCA, wherein the ratio of the thresholds is (FDR to HC-sum):
\begin{equation}
\frac{1 + \sqrt{2 \log \nu p/\widehat{k}}}{\sqrt{2 \rho(\beta)} \sqrt{\log p}}.
\end{equation}
Because of involvement of $\rho(\beta)$, this quantity is hard to analyze. If in an oracle manner, FDR-SEPCA obtained $\widehat{k}$ correctly as $p^{1 - \beta}$, we would find that this ratio is always larger than $1$ for $p > 1$. That is if $\widehat{k}$ assumes the the correct value, HC-sum-SEPCA is  more powerful than FDR-SEPCA. Alternatively, we can note that $\rho(\beta) \in (0, 1]$ and ask when the ratio is larger than $1$. Based on the ratio above, we can see that in the following scenarios
\begin{gather}
\left\lbrace\begin{array}{lcl}
\widehat{k} = 1 & \text{and} & p > 1,\\
2 \leq \widehat{k} \leq 18 & \text{and} & p \geq \widehat{k}\,\text{(always)},\\
\widehat{k} \geq 19 & \text{and} & \log p \geq \frac{1}{8} \left(4 (\log \widehat{k})^2 - 4 \log \widehat{k} + 1\right),
\end{array}\right.
\end{gather}
HC-sum-SEPCA is more powerful than FDR-SEPCA. 

To summarize, we prefer the FDR-controlling alternatives to sum-SEPCA, but depending on the output of FDR-SEPCA, HC-sum-SEPCA may be more powerful. However, as the simulations in Section \ref{sec:sim} revealed (see Figure \ref{fig:risk}), the number of false positives with HC-sum-SEPCA may be higher than with FDR-SEPCA.

\subsection{Overall Message}

We have seen that for $n$ and $p$ sufficiently large and $\vv$ that is sufficiently dense (in the sense of $\|\vv\|_1$ being large), a sum-based statistic and algorithm leads to better performance. This is expected behavior, as by using a sum-based method, we are taking advantage of the equisigned nature of $\vv$. Moreover, within the class of sum-based algorithms, controlling the FDR leads to greater power, as expected. It is difficult to clearly identity which of HC-sum- and FDR-SEPCA will have the greatest power, and the end result may come down to a practitioner's tolerance for false discoveries.

\section{Conclusions}\label{sec:conclusions}

We have considered the setting where the left singular vector of the underlying rank one signal matrix plus noise data matrix is assumed to be sparse and the right singular vector is assumed to be equisigned. We have proposed six different SEPCA algorithms for estimating the sparse principal component based on different decision  statistics and provided sparsistency conditions for the same. Our analysis reveals conditions where a coordinate selection scheme based on a sum-based decision statistic outperforms schemes that utilize the $\ell_1$ and $\ell_2$ decision statistics. Thereby, the proposed algorithm outperforms known schemes such as \textit{diagonal thresholded PCA} \citep{johnstone2009consistency} in terms of estimation of the singular vectors associated with the rank-1 component. We have derived lower bounds on the size of detectable coordinates of the principal left singular vector, utilized these lower bounds to derive lower bounds on the worst-case risk and verified our findings with numerical simulations. Finally, we have discussed the results of our simulations analytically,  by providing a geometric interpretation of the differences in power among the algorithms. 

We note that while we have stated our results for Gaussian noise with identity covariance, we can extend the FWER-controlling results to any log-concave noise distribution, and the FDR-controlling procedures to Gaussian noise with certain non-identity covariances. Additionally, another way to view this work is that it proposes a two-stage procedure/framework for sparse PCA based around hypothesis testing of statistics associated to each row. Some natural extensions would be the inclusion or consideration of other testing frameworks, e.g., that in \cite{ndaoud2018interplay}, where knowledge of the size of coordinates is taken into account. 

\appendices

\section{Proof of Theorem \ref{thm:sparsistency}} \label{sec:thm1_pf}

\renewcommand{\theenumi}{\alph{enumi}}
\begin{enumerate}[leftmargin = *]
\item
Note that 
$$\mathbb{P}\left(T_i \geq \tau\right) \leq \mathbb{P}\left(\max_{j \in I^c} T_j \geq \tau\right)$$
for $i \in I^c$. Taking the maximum over the left-hand side and noting that the right-hand side has limit zero yields the result. This follows from (\ref{eqn:falsepos}). \qed
\item
We consider when true positives occur with probability approaching $1$. We want to find the smallest coordinate $(\theta u_i)$ such that the following probability approaches $1$:
\begin{equation}
\mathbb{P}\left(T_i > \tau_{n, p}\right) = \mathbb{P}\left(\frac{T_i - \mathbb{E} T_i}{\sqrt{\text{Var } T_i}} > \frac{\tau_{n, p} - \mathbb{E} T_i}{\sqrt{\text{Var } T_i}}\right).
\end{equation}
Note that if $(\tau_{n, p} - \mathbb{E} T_i)$ is negative and not tending to zero as $n$ grows, and if the variance of $T_i$ decays to zero as $n$ grows, the quantity
\begin{equation} \label{eqn:tau_inf}
\frac{\tau_{n, p} - \mathbb{E} T_i}{\sqrt{\text{Var } T_i}}
\end{equation}
tends toward negative infinity. Hence, we will specify conditions so that $\text{Var } T_i$ decays to zero as $n$ grows and then compute when a coordinate is detectable by considering when $\tau_{n, p}$ is strictly less than $\mathbb{E} T_i$. For brevity, we omit the computations in solving $\tau_{n, p} < \mathbb{E} T_i$ for $|\theta u_i|$ and present verifications that the variance of $T_i$ has limit $0$. These results show that above the \emph{decision boundary}, we have uniform detection. 

In sum-SEPCA, $T_i$ is a Gaussian random variable with mean $\frac{(\theta u_i)}{\sqrt{n}} \sum_k v_k$ and variance $\frac{\sigma^2}{n}$. Since $\sigma$ does not grow with $n$, $\text{Var } T_i$ always decays to zero. 

In $\ell_2$-SEPCA, $T_i$ has 
$$\mathbb{E} T_i = \left(\theta u_i\right)^2  + \sigma^2$$ 
and 
$$\text{Var } T_i = \frac{2 \sigma^2}{n} \left(\sigma^2 + 2 \left(\theta u_i\right)^2 \right).$$
Since $\sigma$ and $\theta$ are fixed, the variance always decays to $0$.

Let $x_{i,k} = \left(\sqrt{n} \frac{v_k \left(\theta u_i\right)}{\sigma}\right)$. In $\ell_1$-SEPCA, $T_i$ has
\begin{align*}
\text{Var } T_i & = \frac{\sigma^2}{n^2} \sum_k x_{i,k}^2 \left(1 - \left(\text{Erf }\left(\frac{x_{i,k}}{\sqrt{2}}\right)\right)^2\right) \\
& + \frac{\sigma^2}{n} \left(1 - \frac{2}{n \pi} \sum_k \exp\left( -x_{i,k}^2\right) \right) \\
& - 2 \sqrt{\frac{2}{\pi}} \frac{\sigma^2}{n^2} \sum_k x_{i,k}\exp\left(-\frac{x_{i,k}^2}{2}\right) \text{Erf }\left(\frac{x_{i,k}}{\sqrt{2}}\right),
\end{align*}
which is less than or equal to
\begin{equation}
\frac{\left(\theta u_i\right)^2}{n} \sum_k v_k^2  + \frac{\sigma^2}{n} + 2 \sqrt{\frac{2}{\pi}} \frac{\sigma}{n \sqrt{n}} |\left(\theta u_i\right)| \sum_k |v_k|.
\end{equation}

Since $\|\vv\|_2 = 1$, the variance of $T_i$ has limit $0$. Because we cannot solve the inequality $\tau_{n, p} < \mathbb{E} T_i$ analytically, we leave the bound in the form given previously. \qed

In the proof above, note that if $(\tau_{n, p} - \mathbb{E} T_i)$ is positive and not tending to zero as $n$ grows, the quantity in (\ref{eqn:tau_inf}) tends to positive infinity when the variance decays to zero. Hence, modifying the proof by solving $\tau_{n, p} > \mathbb{E} T_i$ for $|\theta u_i|$ yields when a coordinate is not detectable with probability approaching $1$: i.e., when $|\theta u_i|$ is smaller than the values given in (\ref{eqn:limits}). \qed
\end{enumerate}

\subsection{Proof of Corollary \ref{cor:hamming}} \label{ssec:cor_pf}

{The first three terms of (\ref{eqn:expected_hamming}) are characterized by Theorem \ref{thm:sparsistency}(b) and by noting that $I_0 = I$ in the corollary. The last term can be characterized as follows. In particular, for an algorithm with row test statistics $T_i$ and a threshold $\tau$ from Table \ref{tab:details}, we may write 
\begin{equation}
    \sum_{i \notin I} \Prob\left(i \in \widehat{I}\right) = \sum_{i \notin I} \Prob\left(T_i \geq \tau\right).
\end{equation}
}

{For the $\ell_2$- and $\ell_1$-SEPCA algorithms, we may reindex the $T_i$ according to their order statistics $T_{(i)}$, where
$$\left|T_{(1)}\right| \geq \left|T_{(2)}\right| \geq\cdots \geq \left|T_{(p - s)}\right|,$$
and write 
\begin{equation}
    \sum_{i \notin I} \Prob\left(i \in \widehat{I}\right) = \sum_{i = 1}^{p - s} \Prob\left(T_{(i)} \geq \tau\right).
\end{equation}
Note that there are $p - s$ null entries. For the $\ell_2$- and $\ell_1$-SEPCA algorithms, we  have that (as a consequence of \cite[Thm.~3]{Latala2011})
$$\Prob\left(T_{(i)} \geq \tau\right) \leq \left(\frac{1}{e p}\right)^{\sqrt{i}},$$
so that 
\begin{equation} \label{eq:hamming_orderstat_l1l2}
    \sum_{i \notin I} \Prob\left(i \in \widehat{I}\right) = \sum_{i = 1}^{p - s} \Prob\left(T_{(i)} \geq \tau\right) \leq \sum_{i = 1}^{p - s} \left(\frac{1}{e p}\right)^{\sqrt{i}}.
\end{equation}
The right-hand side of (\ref{eq:hamming_orderstat_l1l2}) has limit zero, as needed. 
}

{
For the sum-SEPCA algorithm, from (\ref{eq:thresh_const_2}), it follows that there exists a non-zero constant $\epsilon > 0$ such that the threshold $\tau$ satisfies
$$\tau = \sigma \left(\sqrt{2} + \epsilon\right) \sqrt{\frac{\log p}{n}}.$$
In particular, from (\ref{eq:thresh_const_2}), we have that
$$\epsilon = \frac{1 + \sqrt{\log p} / 3}{\sqrt{2} \textrm{Erf}^{-1}\left(1 - 1/p\right)},$$
where for $p > 1$, 
$$\frac{1}{3 \sqrt{2}} < \epsilon < 2.$$
Hence, for any $T_i$ such that $i \notin I$, 
$$\Prob\left(T_i \geq \tau\right) \leq \exp\left(-\frac{\left(\sqrt{2} + \epsilon\right)^2}{2} \log p\right) = p^{-\left(1 + \sqrt{2} \epsilon + \epsilon^2 / 2\right)},$$
where we have used a Gaussian tail bound \cite[Sec.~2.3]{boucheron2013concentration}. Then,  
\begin{equation}\label{eq:hamming_orderstat_sum}
\begin{split}
    \sum_{i \notin I} \Prob\left(i \in \widehat{I}\right) &= \sum_{i \notin I} \Prob\left(T_i \geq \tau\right) \\ &\leq \left(p - s\right) p^{-\left(1 + \sqrt{2} \epsilon + \epsilon^2 / 2\right)} \\ &\leq p^{-\left(\sqrt{2} \epsilon + \epsilon^2 / 2\right)}.
\end{split}
\end{equation}
Since $\epsilon$ is larger than $1 / 3 \sqrt{2}$, the right-hand side of (\ref{eq:hamming_orderstat_sum}) is upper bounded by $p^{-13/36}$, which has limit zero, as desired.
}

\section{Proof of Theorem \ref{thm:sparsistency_fdr}} \label{sec:thm2_pf}

\subsection{Size of Detectable Coordinates}

\subsubsection{Sum: HC-SEPCA}

If $\vv$ is equisigned, summing across the rows of $\XX$ yields a normally distributed quantity with mean $(\theta u_i) \|\vv\|_1$ and variance $\sigma^2$. Dividing by $\sigma$ and adopting the notation of $HC$, we have that under the alternative hypothesis, $\mu_i = \sqrt{2 r \log p}$, so that 
$$r = \left(\frac{|\theta u_i| \|\vv\|_1}{\sqrt{2 \log p}}\right)^2.$$
Rearranging the inequality $r > \rho(\beta)$ yields
\begin{equation} \label{eqn:HC_sum}
|\theta u_i| > \sigma \sqrt{\rho(\beta)} \frac{\sqrt{2 \log p}}{\|\vv\|_1}.
\end{equation}
Note that sum-SEPCA can detect coordinates of size
\begin{equation} \label{eqn:HC_sum_1}
|\theta u_i| > \sigma C_U \frac{\sqrt{\log p}}{\|\vv\|_1}.
\end{equation}
However, $C_U$ is strictly larger than $\sqrt{2} + 1/(3\sqrt{2})$. Thus, using HC yields a threshold of the same order, but with a strictly smaller scaling. 

\subsubsection{Sum of squares: HC-$\ell_2$-SEPCA}

If we sum the squares of the entries of rows of $\XX$, abusing notation slightly and using $\mathcal{N}(\mu, \sigma^2)$ to indicate a Gaussian random variable with mean $\mu$ and variance $\sigma^2$, the statistic for the $i^{th}$ coordinate is of the form
$$\sum_{k = 1}^n \left(\frac{\sigma}{\sqrt{n}} \mathcal{N}\left(\frac{\theta u_i}{\sigma} v_k \sqrt{n}, 1\right)\right)^2.$$
Assuming oracular knowledge of $\sigma$, the statistic 
$$\frac{n}{\sigma^2} \sum_k X_{ik}^2$$
places us in the setting of (\ref{eqn:HC_chi_test}). The non-centrality parameter $\delta$ is given by 
$$\delta = \sqrt{\sum_{k = 1}^n \left( \frac{\theta u_i}{\sigma} v_k \sqrt{n} \right)^2} = \left|\frac{\theta u_i}{\sigma}\right| \sqrt{n}.$$
Setting $\delta = 2 r \log p$ and solving $r > \rho(\beta)$ yields
\begin{equation} \label{eqn:HC_chi}
|\theta u_i| > \sigma {\rho(\beta)} {\frac{2 \log p}{\sqrt{n}}}.
\end{equation}
We have that $\ell_2$-SEPCA can detect coordinates with
\begin{equation} \label{eqn:HC_chi_1}
|\theta u_i| > \sigma \sqrt{e \sqrt{2}} \sqrt{\frac{1 + \log p}{\sqrt{n}}}.
\end{equation}
Using HC offers a significant improvement over $\ell_2$-SEPCA. However, we also expect HC with the $\chi_n^2$ statistic to have a smaller detectable coordinate: $\|\vv\|_1 \leq \sqrt{n}$, so that for fixed $\beta$ and $p$, the threshold in (\ref{eqn:HC_sum}) is asymptotically larger than that in (\ref{eqn:HC_chi}) (but potentially of the same order). This result is strange in context of the non-FDR results. In any case, HC improves on $\ell_2$-SEPCA. 

\subsubsection{FDR-SEPCA}

Recall that taking sums across the rows of $\XX$, we obtain a vector $\yy$ where $y_i = \mu_i + \sigma z_i$, with $\mu_i = (\theta u_i) \|\vv\|_1$. Moreover, we have noted that 
$$t_k \approx \sqrt{\zeta} (1 + \sqrt{2 \log(\nu  p / k)}),$$ 
where $t_k$ is the level at which $\yy$ is thresholded. It follows that, entries of $\yy$ that are of size at least 
$$y_i > (1 - o(1)) \sqrt{\zeta} \sigma \left(1 + \sqrt{2 \log(\nu  p / \widehat{k})}\right)$$
are selected, or, since $\mu_i = (\theta u_i) \|\vv\|_1$ (when $v$ is equisigned), if we select $\widehat{k}$ coordinates, we expect to detect 
\begin{equation}
\begin{split}
|\theta u_i| &>  (1 - o(1)) \sqrt{\zeta} \sigma \frac{\left(1 + \sqrt{2 \log(\nu  p / \widehat{k})}\right)}{\|v\|_1} \\ &= O\left(\sigma \frac{\sqrt{2 \log \left(\nu  p / \widehat{k}\right)}}{\|\vv\|_1}\right).
\end{split}
\end{equation}
Relative to HC and sum-SEPCA, the gain here is found when there are many smaller coordinates of $\uu$ and $\widehat{k}$ is large. 

\subsection{Proofs for the Higher Criticism-Based Methods}

\renewcommand{\theenumi}{\alph{enumi}}
\begin{enumerate}[leftmargin = *]
\item
From (2.8) in \cite{hall2010innovated}, 
$$\mathbb{P}\left(T_i \geq \tau\right) \leq \mathbb{P}\left(\max_{j \in I^c} T_j \geq \tau\right)$$
has limit zero. \qed
\item
Let $I_1 \subseteq I$ be the set of coordinates with signal larger than the detection limit ($i \in I$ such that $|\theta u_i| > \beta_{crit} (1 + \epsilon)$), and let $I_2 \subseteq I$ contain the rest of the coordinates ($i \in I$ such that $|\theta u_i| < \beta_{crit} (1 - \epsilon)$). By Theorem 1 in \cite{arias2011global}, the asymptotic power for detecting signals below the detection limit is one, and that for signals below the limit is zero. Hence, for $i \in I_1$,
$$\min_{\substack{i \in I~:~  |\theta u_i| > \beta_{crit} (1 + \epsilon)}} \mathbb{P}\left(\text{$i$ selected}\right) \rightarrow 1,$$ 
and for $i \in I_2$,
$$\max_{\substack{i \in I ~:~ |\theta u_i| < \beta_{crit} (1 - \epsilon)}}\mathbb{P}\left(\text{$i$ selected}\right) \rightarrow 0.$$

As with Theorem 1, we omit the computation of $\beta_{crit}$, as it follows from the discussion in Section \ref{ssec:HC}. \qed
\end{enumerate}

\subsection{FDR-SEPCA}

The details of these computations are in Appendix \ref{ssec:fdr_risk_app}, so we will summarize the properties here. 

\renewcommand{\theenumi}{\alph{enumi}}
\begin{enumerate}[leftmargin = *]
\item
The choice of $\nu = 2^{1 / \omega}$ controls the FDR at level $\omega$ \cite{johnstone2014adaptation}. Choosing $\omega = \omega(p) \rightarrow 0$ as $p \rightarrow \infty$ leads to an asymptotic FDR of zero. I.e., for $i \in I^{c}$,
$$\max_{i \in I^c} \qquad \mathbb{P}\left(\text{$i$ selected}\right) \rightarrow 0. \qed$$ 
\item
Noting that the consistency of estimating the mean vector $\bmu = \left(\theta \|v\|_1\right) \uu$ encompasses the estimation of the support of $\uu$, risk bounds for the estimation of $\bmu$ yield the result. To be precise, if the expected risk $\mathbb{E} \|\bmu - \widehat{\bmu}\|_2^2 \leq B$ for some bound $B$, we expect to detect coordinates of size larger than $B$ and to not detect those smaller than $B$. \qed
\end{enumerate}

\section{Risk bounds under $\ell_q$ sparsity} \label{sec:lq_risk}

In this section, we simultaneously generalize our setting to approximate sparsity, and specify the risk lower-bounds. We omit the $\ell_1$-SEPCA algorithm from consideration. 

Let $\uu \in \mathbb{R}^p$ have unit $\ell_2$-norm and belong to an $\ell_q$ ball with radius $C$ for $q \in (0, 2]$. I.e., 
\begin{equation}
\sum_{i = 1}^p |u_i|^q \leq C^q.
\end{equation}
When $q = 0$, we replace $C^q$ with $s$, the level of `hard' sparsity. We have the following result:
\begin{theorem} \label{thm:risk}
Let 
\begin{equation} \label{eqn:loss_appendix}
L(\widehat{\uu}, \uu) = \left\|\uu - \text{sign}(\langle \uu, \widehat{\uu}\rangle) \widehat{\uu}\right\|_2^2
\end{equation}
be the risk of the estimator $\widehat{\uu}$ of $\uu$, where $\uu$ is as specified in (\ref{eqn:model}) and the estimators are the six algorithms that we have previously described. Then, \begin{enumerate}
    \item 
    sum-, HC-sum, and FDR-SEPCA have expected risks lower-bounded by 
        \begin{equation}
        \mathbb{E} L(\widehat{\uu}, \uu) \geq O\left([C^q - 1] \|\vv\|_1^{-(2 - q)}\right).
        \end{equation} 
    \item
    $\ell_2$-SEPCA has a risk lower-bounded by 
        \begin{equation}
        \mathbb{E} L(\widehat{\uu}, \uu) \geq O\left([C^q - 1] n^{-\frac{1}{2}(1 - q/2)} \right).
        \end{equation}
    \item HC-$\ell_2$-SEPCA has a risk lower-bounded by
        \begin{equation}
        \mathbb{E} L(\widehat{\uu}, \uu) \geq O\left([C^q - 1] n^{-(1 - q/2)} \right).
        \end{equation}
\end{enumerate}
\end{theorem}
The rest of this section contains the proof of Theorem \ref{thm:risk}.

\subsection{Proof of Theorem \ref{thm:risk}}

We construct a `worst-case' sparse $\uu$. Note that $C^q \geq 1$ necessarily, and that if $C \geq p^{1 - q/2}$, every unit norm vector is in the $\ell_q$ ball. Hence, we take $C \in [1, p^{1 - q/2})$. 

Let $\theta$ and $\sigma$ be fixed. We want a sparse vector with several coordinates guaranteed to be missed (the probability of not detecting them is asymptotically $1$). For this vector $\uu$ to be sparse and for the loss to not be $1$, set $u_1$ to be $\sqrt{1 - r_n^2}$, where $r_n^2 = o(1)$, and take $u_2, \cdots, u_{m_n + 1}$ to be ${r_n}/{\sqrt{m_n}}$. The other coordinates of $u$ are $0$, so that $\uu$ has unit $\ell_2$-norm. 

We assume that $u_1$ is detected with probability $1$ as $n \rightarrow \infty$, and want to set $u_2, \cdots, u_{m_n + 1}$ so that the expected loss is lower bounded by:
\begin{equation}
\begin{split}
\mathbb{E} L(\uu, \widehat{\uu}) &\geq \sum_{k = 1}^p |u_k|^2 \mathbb{P}(\text{Not Selecting Coordinate k}) \\
&\geq \sum_{k = 2}^{m_n + 1} |u_k|^2 \mathbb{P}(\text{Not Selecting Coordinate k}).
\end{split}
\end{equation}
If coordinates of size $\frac{r_n}{\sqrt{m_n}}$ are not detected, the expected loss is lower bounded by $r_n^2$. 

Let $m_n = \lfloor m \rfloor$ where 
$$m = \delta n^{\phi} r_n^{\psi} \|v\|^{\eta}.$$
Note that we have not specified the norm used in $\|v\|$: we will choose the norm at the very end of the calculation. Let 
$$r_n = [C^q - 1]^{\alpha} n^{\beta + \gamma q} \|v\|^{\kappa},$$
so that, 
$$\frac{r_n}{\sqrt{m_n}} \approx \frac{r_n}{\sqrt{m}} = \frac{1}{\delta} n^{-\phi / 2} \|v\|^{\kappa - \eta/2} r_n^{1 - \psi/2}.$$
We will choose $\delta, \phi, \eta, \alpha, \beta, \gamma, \kappa, \psi$ so that the $\ell_q$ sparsity constraint is met and the lower bound $r_n^2$ is maximized. The sparsity constraint requires: 
\begin{gather}
\begin{split}
\sum_{i = 1}^p |u_i|^q = (1 - r_n)^{q/2} + m_n^{1 - q/2} r_n^q  \leq 1 + m^{1 - q/2} r_n^q  \leq C^q.
 \end{split}
\end{gather}

First, we will assume (for now) that $r_n = o(1)$ and that via other parameters we may control the scaling of the coordinate sizes; hence, we set $\psi = 2$. Then, 
\begin{equation}
\begin{split}
r_n^q m^{1 - q/2} = \delta^{1 - q/2} n^{2(\beta + \gamma q) + (1 - q/2) \phi} [C^q - 1]^{2 \alpha} \|v\|^{2 \kappa + \eta(1 - q/2)}.
\end{split}
\end{equation}
We need this quantity to be smaller than $C^q - 1$. To eliminate the $n$ dependence, we set $\beta = \frac{-\phi}{2}$ and $\gamma = \frac{\phi}{4}$. We choose $\alpha = \frac{1}{2}$ to match powers of $[C^q - 1]$ on both sides of the inequality.  Defining another parameter $\rho$, let $\delta = \rho \|\vv\|^{-\eta}$. Then, the inequality is 
$$\rho^{1 - q/2} \|\vv\|^{2 \kappa} [C^q - 1] \leq [C^q - 1].$$ 
Choosing $\rho \leq \|\vv\|^{-2 \kappa / (1 - q /2)}$ is enough. 

With these choices of parameters, 
$$r_n = \sqrt{[C^q - 1]} n^{-\frac{1}{2} \phi (1 - q/2)} \|\vv\|^{\kappa},$$ 
and 
$$m = \rho n^{\phi} r_n^2,$$
so that 
$$\frac{r_n}{\sqrt{m}} = \frac{1}{\sqrt{\rho}} n^{-\phi /2}.$$
Noting that 
$$\frac{1}{\sqrt{\rho}} \geq \|\vv\|^{\kappa / (1 - q/2)},$$
choosing $\rho = \|\vv\|^{-2 \kappa / (1 - q /2)}$ leads to the smallest possible choice of coordinate. 

In summary:
\begin{equation}
r_n = \sqrt{[C^q - 1]} n^{-\frac{1}{2} \phi (1 - q/2)} \|\vv\|^{\kappa},
\end{equation}
\begin{equation}
r_n^2 = [{C^q - 1}] n^{-\phi (1 - q/2)} \|\vv\|^{2 \kappa},
\end{equation}
\begin{equation}
m =  \|\vv\|^{-2 \kappa / (1 - q /2)} n^{\phi} r_n^2,
\end{equation}
and
\begin{equation}
\frac{r_n}{\sqrt{m}} = \|\vv\|^{\kappa / (1 - q /2)} n^{-\phi /2}.
\end{equation}

So, for a given algorithm, it remains to choose $\phi$ and $\kappa$ so that the worst-case risk is lower-bounded by $r_n^2$. Sum-SEPCA misses coordinates of size $O\left(\frac{\sqrt{\log p}}{\|\vv\|_1}\right)$ and $\ell_2$-SEPCA misses coordinates of size $O\left(\frac{\sqrt{\log p}}{n^{1/4} \|\vv\|_2}\right)$. For sum-SEPCA, $\kappa = \frac{q - 2}{2}$, and for $\ell_2$-SEPCA, $\kappa$ is irrelevant, as $\|\vv\|_2 = 1$. Sum-SEPCA uses $\phi = 0$ and $\ell_2$-SEPCA uses $\phi = \frac{1}{2}$. Hence, sum-SEPCA has a risk lower-bounded by 
\begin{equation}
O\left([C^q - 1] \|\vv\|_1^{-(2 - q)}\right).
\end{equation} 
Noting that $\|\vv\|_2 = 1$, $\ell_2$-SEPCA has 
\begin{equation}
\begin{split}
& O\left([C^q - 1] n^{-\frac{1}{2}(1 - q/2)} \|\vv\|_2^{-(1 - q/2)}\right) = \\& O\left([C^q - 1] n^{-\frac{1}{2}(1 - q/2)} \right).
\end{split}
\end{equation}

In the $\ell_0$ case, i.e., when $\uu$ has no more than $s$ non-zero entries, the preceding analysis goes through with $C^q$ replaced by $s$ and $q$ set to zero. 

\subsubsection{FDR Algorithms}

For HC-sum-SEPCA, the $\beta_{crit}$ is of the same order as that for sum-SEPCA. Similarly, for FDR-SEPCA, if $\widehat{k}$ is much smaller than $p$, $\beta_{crit}$ is of roughly the same order. Hence, these two algorithms have the same risk bound as sum-SEPCA. For HC-$\ell_2$-SEPCA, $\kappa = 0$ and $\phi = 1$. The risk is therefore lower-bounded by 
\begin{equation}
O\left([C^q - 1] n^{-(1 - q/2)} \right).
\end{equation}

\section{FDR-SEPCA: Further Details} \label{sec:fdr_details}

Let $y_i = \mu_i + \sigma z_i$, where $i \in \lbrace 1, \cdots, p\rbrace$, the vector $\zz$ of the $z_i$ is normally distributed with mean $0$ and covariance $\bSig$, and $\bSig$ satisfies 
$$\xi_o \mathcal{I}_p \leq \bSig \leq \xi_1 \mathcal{I}_p.$$
Here, $\xi_0$ is the smallest eigenvalue of $\bSig$ and $\xi_1$ is the largest. The mean vector $\bmu$ of the $\mu_i$ is assumed to be sparse; the goal is to estimate $\bmu$. The following penalized least squares formulation yields an estimator for $\bmu$: 
\begin{equation} \label{eqn:FDR_prob}
\widehat{\bmu} = \arg \min_{\bmu} \|\yy - \bmu\|_2^2 + \sigma^2 \text{pen}\left(\|\bmu\|_0\right),
\end{equation}
where $\text{pen}(k)$ is defined as 
\begin{equation}
\text{pen}(k) = \xi_1 \zeta k \left(1 + \sqrt{2 L_{p, k}}\right)^2,
\end{equation}
with $\zeta > 1$ and 
\begin{equation} 
L_{p, k} = (1 + 2 \beta) \log(\nu p / k).
\end{equation}
The parameter $\beta$ may be set to $0$ here, and $\nu$ is chosen to be no smaller than $e^{1/(1 + 2 \beta)}$. We define $\|\bmu\|_0$ to be number of non-zero coordinates of $\bmu$. 

The solution to (\ref{eqn:FDR_prob}) is given by hard-thresholding. Let $|y|_{(i)}$ be the $i^{th}$ order statistic of $|y_i|$, namely $|y|_{(1)} \geq \cdots \geq |y|_{(p)}$. Then if
\begin{equation}
\widehat{k} = \arg \min_{k \geq 0} \sum_{i > k} |y|_{(i)}^2 + \sigma^2 \text{pen}(k), 
\end{equation}
defining 
\begin{equation}
t_k^2 = \text{pen}(k) - \text{pen}(k - 1),
\end{equation}
the solution is to hard threshold at $t_{\widehat{k}}$. 

In this set-up, we have that 
$$t_k \approx \lambda_{p, k} = \sqrt{\xi_1 \zeta} (1 + \sqrt{2 L_{p,k}}),$$
with $|t_k - \lambda_{p, k}| < c / \lambda_k$. More precisely, Lemma 11.7 of \cite{johnstone2013gaussian} says that 
$$\lambda_{p,k} - \frac{4 \zeta b}{\lambda_{p,k}} \leq t_k \leq \lambda_{p,k}.$$
When $\nu \geq e^2$, we may take $b = (1 + 2 \beta)$. In any case, if $k = o(n)$, $\lambda_{p, k} \asymp \sqrt{\log p}$. Hence, entries of $\yy$ that are of size at least 
$$y_i > (1 - o(1)) \sqrt{\xi_1 \zeta} \sigma \left(1 + \sqrt{2 \log(\nu p / \widehat{k})}\right)$$
are selected, or, 
since $\mu_i = (\theta u_i) \|\vv\|_1$ (when $\vv$ is equisigned), if we select $\widehat{k}$ coordinates, we expect to detect 
\begin{equation}
\begin{split}
|\theta u_i| &>  (1 - o(1)) \sqrt{\xi_1 \zeta} \sigma \frac{\left(1 + \sqrt{2 \log(\nu p / \widehat{k})}\right)}{\|\vv\|_1} \\ &= O\left(\sigma \sqrt{\xi_1} \frac{\sqrt{2 \log \left(\nu p / \widehat{k}\right)}}{\|\vv\|_1}\right).
\end{split}
\end{equation}

\subsection{Risk Behavior} \label{ssec:fdr_risk_app}

Recalling that (\ref{eqn:FDR_prob}) solves a penalized least squares problem for $\widehat{\yy}$ close to $\yy$, we may discuss the statistical behavior of this estimator. The following discussion follows and reproduces that in \cite{johnstone2014adaptation}

First, note that for $\beta = 0$, the parameter $\nu$ directly controls the FDR (where a false positive corresponds to selecting a zero coordinate in $\yy$): a choice of $\nu = 2^{1/\omega}$ for $\omega \in (0, 1)$ bounds the FDR at a level $\omega$. 

Second, the expected risk, $\mathbb{E} \|\yy - \widehat{\yy}\|_2^2$, is bounded as follows. By Proposition 4.1 in \cite{johnstone2014adaptation}, 
\begin{equation} \label{eqn:expected_risk_fdr}
\mathbb{E} \|\yy - \widehat{\yy}\|_2^2 \leq D\left[2 M_p' \xi_1 \sigma^2 + \mathcal{R}(\yy, \sigma)\right],
\end{equation}
where $D$ is a constant $2 \zeta (\zeta + 1)^3 (\zeta - 1)^{-3} = \Theta(1)$, we assume that $\xi_1 = 1$, and $0 \leq M_p' \leq C_{\beta} p^{-2 \beta} \nu^{-1}$, for some $C_{\beta} > 0$. Since $\beta = 0$, $M_p' = O(1/\nu) = O(\omega)$, if we control the FDR at level $\omega$. 

The second term in (\ref{eqn:expected_risk_fdr}) is the ideal risk, or, the infimum of the penalized least squares objective. If $\yy$ belongs to an $\ell_q$ ball with radius $C$ and $0 < q < 2$, and we define
\begin{equation}
r_{p, q}(C) = \left\lbrace \begin{array}{ll} C^2 & \text{ if } C \leq \sqrt{1 + \log p}, \\ C^q [1 + \log(p / C^q)]^{1 - q/2} & \text{ if } \sqrt{1 + \log p} \leq C \leq p^{1/q}, \\
p & \text{ if } C \geq p^{1/q},\end{array} \right.
\end{equation}
the ideal risk is bounded as
\begin{equation}
\sup_{\yy \in \mathbb{R}^p : \sum_i |y_i|^q \leq C^q} \mathcal{R}(\yy, \sigma) \leq c (\log \nu) \sigma^2 r_{p, q}(C / \sigma),
\end{equation}
for some $c > 0$. The supplementary results in \cite{johnstone2014adaptation} yield that $\mathcal{R}(\yy, \sigma)$ is bounded by $C^2 \log \nu$, and by $C^2$ when $C \leq \sqrt{1 + \log p}$.

As in Appendix \ref{sec:lq_risk}, we may replace $q$ with $0$ and $C^q$ with $s$ in the case of hard sparsity with $s$ non-zero coordinates. Doing so leads to the bound:
\begin{equation}
\mathbb{E}\|\yy - \widehat{\yy}\|_2^2 \leq s \sigma \log \nu \log \frac{\sigma p \nu}{s} + \sigma_1 \sigma^2 \frac{2}{\nu}.
\end{equation}
Note that we have recovered the factor of $\log \nu p / s$ in $\beta_{crit}$. 


\bibliographystyle{plain} 
\bibliography{sepca}

\end{document}